\input amstex
\documentstyle{amsppt}
\magnification1200

\expandafter\ifx\csname bethe.def\endcsname\relax \else\endinput\fi
\expandafter\edef\csname bethe.def\endcsname{%
 \catcode`\noexpand\@=\the\catcode`\@\space}
\catcode`\@=11

\mathsurround 1.6pt

\def\hcor#1{\advance\hoffset by #1}
\def\vcor#1{\advance\voffset by #1}
\let\bls\baselineskip  \let\ignore\ignorespaces
\def\vsk#1>{\vskip#1\bls} \let\adv\advance 
\def\vv#1>{\vadjust{\vsk#1>}\ignore} \def\vvv#1>{\vadjust{\vskip#1}\ignore}
\def\vvn#1>{\vadjust{\nobreak\vsk#1>\nobreak}\ignore}
\def\vvvn#1>{\vadjust{\nobreak\vskip#1\nobreak}\ignore}
\def\setnormalbls{\edef\normalbls{\bls\the\bls}}
\def\setmaths{\edef\maths{\mathsurround\the\mathsurround}}

 \let\nt\noindent 
\def\nn#1>{\noalign{\vskip #1pt}} \def\NN#1>{\openup#1pt}
 
\let\Sum\sum \def\sum{\Sum\limits} 
\let\Prod\prod \def\prod{\Prod\limits} \let\Int\int \def\int{\Int\limits}

\let\=\m@th \def\&{.\kern.1em} \def\>{\!\;} \def\:{\!\!\;}

\ifx\plainfootnote\undefined \let\plainfootnote\footnote \fi
\expandafter\ifx\csname amsppt.sty\endcsname\relax
 
\else \fi

\newbox\s@ctb@x
\def\s@ct#1 #2\par{\removelastskip\vsk>
 \vtop{\bf\setbox\s@ctb@x\hbox{#1} \parindent\wd\s@ctb@x
 \ifdim\parindent>0pt\adv\parindent.5em\fi\item{#1}#2\strut}%
 \nointerlineskip\nobreak\vtop{\strut}\nobreak\vsk-.4>\nobreak}

\newbox\t@stb@x
\def\gadv{\global\advance} \def\gad#1{\gadv#1 1} 
\def\l@b@l#1#2{\def\n@@{\csname #2no\endcsname}%
 \if *#1\gad\n@@ \expandafter\xdef\csname @#1@#2@\endcsname{\the\Sno.\the\n@@}%
 \else\expandafter\ifx\csname @#1@#2@\endcsname\relax\gad\n@@
 \expandafter\xdef\csname @#1@#2@\endcsname{\the\Sno.\the\n@@}\fi\fi}
\def\l@bel#1#2{\l@b@l{#1}{#2}\?#1@#2?}
\def\?#1?{\csname @#1@\endcsname}
\def\[#1]{\def\n@xt@{\ifx\t@st *\def\n@xt####1{{\setbox\t@stb@x\hbox{\?#1@F?}%
 \ifnum\wd\t@stb@x=0 {\bf???}\else\?#1@F?\fi}}\else
 \def\n@xt{{\setbox\t@stb@x\hbox{\?#1@L?}\ifnum\wd\t@stb@x=0 {\bf???}\else
 \?#1@L?\fi}}\fi\n@xt}\futurelet\t@st\n@xt@}
\def\(#1){{\rm\setbox\t@stb@x\hbox{\?#1@F?}\ifnum\wd\t@stb@x=0 ({\bf???})\else
 (\?#1@F?)\fi}}
\def\dff{\expandafter\d@f} \def\d@f{\expandafter\def}
\def\edff{\expandafter\ed@f} \def\ed@f{\expandafter\edef}

\newcount\Sno \newcount\Lno \newcount\Fno
\def\Section#1{\gad\Sno\Fno=0\Lno=0\s@ct{\the\Sno.} {#1}\par} \let\Sect\Section
\def\section#1{\gad\Sno\Fno=0\Lno=0\s@ct{} {#1}\par} \let\sect\section
\def\l@F#1{\l@bel{#1}F} \def\<#1>{\l@b@l{#1}F} \def\l@L#1{\l@bel{#1}L}
\def\Tag#1{\tag\l@F{#1}} \def\Tagg#1{\tag"\llap{\rm(\l@F{#1})}"}
\def\Th#1{Theorem \l@L{#1}} \def\Lm#1{Lemma \l@L{#1}}
\def\Prop#1{Proposition \l@L{#1}}
\def\Cr#1{Corollary \l@L{#1}} \def\Cj#1{Conjecture \l@L{#1}}
 
\def\Proof#1.{\demo{\it Proof #1}}

\def\Par{\par\medskip} \def\setparindent{\edef\Parindent{\the\parindent}}
\def\Appendix{\Sno=64\let\p@r@\z@ 
 \def\Section##1{\gad\Sno\Fno=0\Lno=0 \s@ct{} \hskip\p@r@ Appendix
\char\the\Sno
  \if *##1\relax\else {.\enspace##1}\fi\par} \let\Sect\Section
 \def\section##1{\gad\Sno\Fno=0\Lno=0 \s@ct{} \hskip\p@r@ Appendix%
  \if *##1\relax\else {.\enspace##1}\fi\par} \let\sect\section
 \def\l@b@l##1##2{\def\n@@{\csname ##2no\endcsname}%
 \if *##1\gad\n@@
 \expandafter\xdef\csname @##1@##2@\endcsname{\char\the\Sno.\the\n@@}%
 \else\expandafter\ifx\csname @##1@##2@\endcsname\relax\gad\n@@
 \expandafter\xdef\csname @##1@##2@\endcsname{\char\the\Sno.\the\n@@}\fi\fi}}

  \let\8\infty

\let\=\m@th  \def\_#1{_{\rlap{$\ssize#1$}}}

\def\lc{\lsym,}

\def\E(#1){\mathop{\hbox{\rm End}\,}(#1)} 
\def\id{\hbox{\rm id}}  
\def\tr{\hbox{\rm tr}}

\def\1{^{-1}} \def\vst#1{{\lower2.1pt\hbox{$\bigr|_{#1}$}}}

 \let\eps\varepsilon \let\epsilon\eps
\let\ka\kappa
\let\la\lambda 
  
 \let\phi\varphi

\let\logo@\relax
\let\m@k@h@@d\makeheadline \let\m@k@f@@t\makefootline
\def\makeheadline{\ifnum\pageno=1\headline={\hfil}\fi\m@k@h@@d}
\def\makefootline{\ifnum\pageno=1\footline={\hfil}\fi\m@k@f@@t}



\hsize12.7cm
\vsize19.8cm

\NoRunningHeads

\def\aN{{\frak{a}_N}}

\def\AaNZ{{\operatorname{B}(\glN,\sigma,Z)}}
\def\ANZ{{\operatorname{B}(\glN, Z)}}
\def\AgNZ{{\operatorname{J\hskip0.5pt}(\glN, Z)}}
\def\AgaNZ{{\operatorname{J\hskip0.5pt}(\glN,\sigma,Z)}}

\def\bN{{\frak{b}_N}}

\def\CC{{\Bbb C}}

\def\de{\delta}
\def\deg{{\operatorname{deg}\ts}}
\def\det{{\operatorname{det}\ts}}
\def\dim{{\operatorname{dim}\ts}}

\def\End{{\operatorname{End}}}
\def\EndCN{{\End(\CC^N)}}
\def\EndLak
{{\End\bigl({\Lambda\hskip-2pt}^{\raise0.5pt\hbox{$\scriptstyle k$}}
\hskip-1pt\Bbb C^N\bigr)}}
\def\EndLal
{{\End\bigl({\Lambda\hskip-2pt}^{\raise0.5pt\hbox{$\scriptstyle l$}}
\hskip-1pt\Bbb C^N\bigr)}}
\def\enddemos{{$\quad\square$\enddemo}}
\def\ep{\epsilon}

\def\Fk{{\operatorname{F\hskip-3pt}_{\lower0.8pt\hbox{$\scriptstyle k$}}
\hskip-1.5pt(\Bbb C^N)}}
\def\Fl{{\operatorname{F\hskip-2.5pt}_{\lower0.8pt\hbox{$\scriptstyle l$}}
\hskip-1.5pt(\Bbb C^N)}}

\def\glN{{\frak{gl}_N}}

\def\gp{h}

\def\id{{\operatorname{id}\ts}}
\def\io{\iota}

\def\ka{\kappa}
\def\kN{{\frak{t}_{M,N}}}
\def\kaN{{\frak{s}_{M,N}}}

\def\la{\lambda}
\def\lc{{,\ts\ldots,\ts}}

\def\max{{\operatorname{max}\hskip1pt}}
\def\min{{\operatorname{min}\hskip1pt}}

\def\ot{\otimes}

\def\RR{{\Bbb R}}

\def\sgn{{\operatorname{sgn}}}
\def\so{{\frak{so}}}
\def\soN{{\so_N}}

\def\sp{{\frak{sp}}}
\def\spN{{\sp_N}}

\def\tr{{\operatorname{tr}}}
\def\tR{{\tilde R}}
\def\ts{\thinspace}
\def\tS{{\tilde S}}
\def\tT{{\widetilde T}}

\def\UaN{{\operatorname{U}(\aN)}}
\def\UN{{\operatorname{U}(\glN)}}
\def\uMN
{{\frak{g\hskip-2pt}_{\lower0.8pt\hbox{$\scriptstyle M,N$}}}}
\def\uMNs
{{\frak{g\hskip-2pt}_{\lower0.8pt\hbox{$\scriptstyle M,N$}}^{\hskip3pt\ast}}}
\def\uaMN
{{\frak{f\hskip-2pt}_{\lower0.8pt\hbox{$\scriptstyle M,N$}}}}
\def\uaMNs
{{\frak{f\hskip-2pt}_{\lower0.8pt\hbox{$\scriptstyle M,N$}}^{\hskip3pt\ast}}}
\def\UNx{{\operatorname{U}(\glN[\hskip1pt t \hskip1pt])}}
\def\UaNx{{\operatorname{U}(\aN[\hskip1pt t \hskip1pt])}}

\def\wA{{\widehat A}}
\def\wB{{\widehat B}}
\def\wS{{\widehat S}}
\def\wT{{\widehat T}}

\def\xN{{\frak{x}_N}}

\def\YN{{\operatorname{Y}(\glN)}}
\def\YgN{{\operatorname{X}(\glN)}}
\def\YgaN{{\operatorname{X}(\glN,\sigma)}}
\def\YgoN{{\operatorname{X\hskip-1pt}_{\lower0.6pt\hbox{$\scriptstyle 1$}}
\hskip-1pt(\glN)}}
\def\YgtN{{\operatorname{X\hskip-1pt}_{\lower0.6pt\hbox{$\scriptstyle 2$}}
\hskip-1pt(\glN)}}
\def\YgtaN{{\operatorname{X\hskip-1pt}_{\lower0.6pt\hbox{$\scriptstyle 2$}}
\hskip-1pt(\glN,\sigma)}}
\def\YgnN{{\operatorname{X\hskip-1pt}_{\lower0.6pt\hbox{$\scriptstyle M$}}
\hskip-1pt(\glN)}}
\def\YganN{{\operatorname{X\hskip-1pt}_{\lower0.6pt\hbox{$\scriptstyle M$}}
\hskip-1pt(\glN,\sigma)}}
\def\YgnoN{{\operatorname{X\hskip-1pt}_{\lower0.6pt\hbox{$\scriptstyle M+1$}}
\hskip-1pt(\glN)}}
\def\YganoN{{\operatorname{X\hskip-1pt}_{\lower0.6pt\hbox{$\scriptstyle M+1$}}
\hskip-1pt(\glN,\sigma)}}
\def\YrN{{\operatorname{Y\!}_r(\glN)}}
\def\YrmN{{\operatorname{Y\!}_{r-1}(\glN)}}
\def\YpqN{{\operatorname{Y\!}_{p+q-2}(\glN)}}
\def\YaN{{\operatorname{Y}(\glN,\sigma)}}

\def\ZZ{{\Bbb Z}}

\font\bigbf=cmbx10 scaled 1200

\line{\bigbf Bethe Subalgebras in Twisted Yangians\hfill}
\bigskip\bigskip
\line{{\bf Maxim Nazarov\hskip0.2pt}$^1$
and \,\,{\bf Grigori Olshanski\hskip1pt}$^2$\hfill}
\bigskip
\line{${}^1$ University of Wales, Swansea SA2 8PP, United Kingdom\hfill}
\line{\phantom{${}^1$} E-mail address: mamln\@ swansea.ac.uk\hfill}
\line{${}^2$ Institute for Problems of Information Transmission,
Moscow 101447, Russia\hfill}
\line{\phantom{${}^2$} E-mail address: olsh\@ ippi.ac.msk.su\hfill}

\bigskip\bigskip
\noindent
{\bf Abstract:}
We study analogues of the Yangian of the Lie algebra $\frak{gl}_N$
for the other classical Lie algebras $\frak{so}_N$ and $\frak{sp}_N$.
We call them twisted Yangians. They are {coideal} subalgebras
in the Yangian of $\frak{gl}_N$
and admit
homomorphisms onto the universal enveloping algebras
$\operatorname{U}(\frak{so}_N)$ and $\operatorname{U}(\frak{sp}_N)$
respectively. In every twisted Yangian we
construct a family of maximal commutative subalgebras
parametrized by the regular semisimple elements of
the corresponding classical Lie algebra. The images 
in $\operatorname{U}(\frak{so}_N)$ and $\operatorname{U}(\frak{sp}_N)$
of these subalgebras are also maximal commutative.

\bigskip\bigskip
\line{\bf Introduction\hfill}
\medskip
\nt
In this article we study the Yangian of the Lie algebra $\glN$
and its analogues for the other classical Lie algebras $\soN$ and $\spN$.
The Yangian $\YN$ is a deformation of the
universal enveloping algebra $\UNx$ in the class of Hopf algebras [D1].
Moreover, it contains the universal enveloping algebra
$\UN$ as a subalgebra and admits a homomorphism
$\pi:\YN\to\UN$ identical on $\UN$.

Let $\aN$ be one of the Lie algebras $\soN$ and $\spN$.
In [D1] the Yangian $\operatorname{Y}(\aN)$ was defined as a deformation
of the Hopf algebra $\UaNx$. It contains $\UaN$ as a subalgebra
but does not admit a homomorphism
$\operatorname{Y}(\aN)\to\UaN$ identical on $\UaN$.
In the present article we consider another analogue of the Yangian $\YN$
for the classical Lie algebra $\aN$. It has been introduced
in [O2] and called twisted Yangian; see also 
[MNO]. The definition in [O2]
was motivated by [O1] and [C2\ts,\ts S]. Algebras closely related to this
analogue of $\YN$ were recently studied in [NS].

Consider $\aN$ as a fixed point subalgebra in the Lie algebra $\glN$
with respect to an involutive automorphism $\sigma$. The twisted Yangian
$\YaN$ is a subalgebra in $\YN$. Moreover, it is a left coideal in
the Hopf algebra $\YN$. It also contains $\UaN$ as a subalgebra
and does admit a homomorphism $\rho:\YaN\to\UaN$ identical on $\UaN$;
see Section 3.
The algebra $\YaN$ is a deformation of the universal enveloping
algebra for the twisted current Lie algebra
$$
\bigl\{\ts F(t)\in\glN[\hskip1pt t \hskip1pt]
\ts\ts|\ts\ts
\sigma\bigl(F(t)\bigr)=F(-t)\ts\bigr\}.
$$

There is a remarkable family of maximal commutative subalgebras
in $\YN$.
They are parametrized by the regular semisimple
elements of $\glN$. As well as the Yangian $\YN$ itself, these subalgebras
were studied in the works
by mathematical physicists from St.\ts Petersburg
on Bethe Ansatz;
see for instance [KBI] and [KR\ts,\ts KS].
These subalgebras were also studied in [C1].
We will call them
Bethe subalgebras. In Section~1 of the present article we recall their
definition. Their images in $\UN$
with respect to the homomophism $\pi$ are also maximal commutative;
see Section 2.

\newpage

The main aim of this article is to construct analogues of
the Bethe subalgebras in $\YN$ for the twisted Yangian $\YaN$.
In Section 3 for any element $Z\in\aN$ we construct certain commutative
subalgebra
in $\YaN$. This construction is a generalization
of one result from [S]. If the element $Z\in\aN$
is regular semisimple then the corresponding commutative subalgebra
in $\YaN$ is maximal.
\text{Moreover,} the image of this subalgebra in $\UaN$
with respect to the homomorphism $\rho$ is also maximal commutative;
see Section 4. This image
in $\UaN$
is a quantization 
of a maximal involutive subalgebra in the Poisson
algebra $\operatorname{S}(\aN)$ obtained by the so called
shift of argument method;
see [K2] and [MF]. 
For further details on the involutive subalgebras
in $\operatorname{S}(\aN)$ obtained by this method
see for instance [RS]. Some results on the quantization of these
subalgebras
can be found in [V].
\Par
We are indebted to M.\,Ra{\"\i}s who explained to us that
the methods of [K1] can be applied to the current Lie algebras; see [RT].
Together with [M] and [MNO] the present article is a part of
a project on representation theory of Yangians initiated by [O1\ts,\ts O2].
It is our joint project
with A.\,Molev, and we are grateful to him for
collaboration.

\Section{Bethe subalgebras in Yangians}
\nt
We will start this section with recalling several known facts
from [D1] and [KR\ts,\ts KS] about
the {\it Yangian} of the Lie algebra $\glN$; see also [MNO\ts,\ts Sections 1--2].
This is a complex associative
unital algebra $\YN$ with the countable
set of generators $T^{(r)}_{ij}$ where $r=1,2,\ts\ldots$ and
$i,j=1,\ts\dots\ts,N$. The de\-fin\-ing relations in the algebra $\YN$ are
$$
[\ts T_{ij}^{(p+1)},T_{kl}^{(q)}\ts]-
[\ts T_{ij}^{(p)},T_{kl}^{(q+1)}\ts]=
T_{kj}^{(p)}\ts T_{il}^{(q)}-T_{kj}^{(q)}\ts T_{il}^{(p)};
\qquad
p,q=0,1,2,\ldots\ts
\Tag{1.1}
$$
where $T_{ij}^{(0)}=\de_{ij}\cdot1$.
The collection \(1.1) is equivalent to the collection of relations
$$
[\ts T_{ij}^{(p)},T_{kl}^{(q)}\ts]=
\sum_{r=1}^{\min(p,q)}
\bigl(\ts
T_{kj}^{(r-1)}\ts T_{il}^{(p+q-r)}-
T_{kj}^{(p+q-r)}\ts T_{il}^{(r-1)}
\ts\bigr);
\quad
p,q=1,2,\ldots\ts.
\Tag{1.11}
$$

Let $E_{ij}\in\EndCN$ be the standard matrix units.
We will also use the following matrix form of the relations \(1.1).
Introduce a formal variable $u$ and consider the
{\it Yang $R$-matrix}
$$
R(u)=u\cdot\id-\sum_{i,j}\ts E_{ij}\ot E_{ji}
\in\EndCN^{\ot\ts2}[u]\ts,
$$
the indices $i$ and $j$ run through the set $\{1\lc N\}$.
We will employ the equality
$$
R(u)\ts R(-u)=(\ts1-u^2)\cdot\id.
\Tag{1.0}
$$
Introduce the formal power series in $u\1$
$$
T_{ij}(u)=
T_{ij}^{(0)}+T_{ij}^{(1)}\ts u\1+T_{ij}^{(2)}\ts u^{-2}+\ldots
$$
and combine all these series into the single element
$$
T(u)=\sum_{i,j}\ts E_{ij}\ot T_{ij}(u)
\in\EndCN\ot\YN\ts[[u\1]].
$$

We will also regard $E_{ij}$ as generators of the universal
enveloping algebra $\UN$.
The algebra $\YN$ contains $\UN$ as
a subalgebra: due to \(1.1) the assignment $E_{ij}\mapsto T_{ij}^{(1)}$
defines the embedding. Moreover, there is a homomorphism
$$
\pi:\ts\YN\to\UN:\ts T_{ij}(u)\mapsto \de_{ij}+E_{ij}\ts u^{-1}.
\Tag{1.1005}
$$
The homomorphism $\pi$ is by definition
identical on the subalgebra $\UN$. It is
\line{called
the {\it evaluation homomorphism} for the algebra $\YN$.
There is a natural} Hopf algebra structure on $\YN$. The comultiplication
$\Delta:\YN\to\YN^{\ot2}$ is defined by the assignement
$$
T_{ij}(u)\mapsto\sum_{k}\ts
T_{ik}(u)\ot T_{kj}(u).
\Tag{1.*}
$$
Here the tensor product is taken over the subalgebra $\CC[[u^{-1}]]$
in $\YN\ts[[u\1]]$.

Throughout this article we will denote by by $\io_s$ the
embedding of the algebra $\EndCN$ into a finite tensor product
$\EndCN^{\ot\ts n}$
as the $s$-th tensor factor:
$$
\io_s(X)=1^{\ot\ts (s-1)}\ot X\ot1^{\ot\ts(n-s)};\qquad s=1,\dots,n.
$$
For any series $Y(u)\in\EndCN\ot\YN[[u\1]]$ we set
$$
Y_s(u)=\io_s\ot\id\bigl(Y(u)\bigr)\in\EndCN^{\ot n}\ot\YN[[u\1]].
$$
Let $v$ be another formal variable.
In the above notation the defining relations \(1.1) can be rewritten as
the single relation in
$\EndCN^{\ot2}\!\ot\YN\ts((u\1,v\1))$
$$
R(u-v)\ot1\cdot\ts T_1(u)\ts T_2(v)=
T_2(v)\ts T_1(u)\cdot R(u-v)\ot1.
\Tag{1.2}
$$

Let $\tr_n$ be the standard matrix trace on $\EndCN^{\ot n}$.
We will also use various embeddings
of the algebra $\EndCN^{\ot\ts m}$ into $\EndCN^{\ot\ts n}$
for any $m\leqslant n$.
When $1\leqslant s_1<\ldots<s_m\leqslant n$ and $X\in\EndCN^{\ot m}$
we put
$$
X_{s_1\ldots s_m}=
\io_{s_1}\!\ot\ldots\ot\io_{s_m}(X)\in\EndCN^{\ot n}.
$$
For any series $Y(u)\in\EndCN^{\ot\ts m}\ot\YN[[u\1]]$ we will also put
$$
Y_{s_1\ldots s_m}(u)=
\io_{s_1}\!\ot\ldots\ot\io_{s_m}\ot\id\bigl(Y(u)\bigr)
\in\EndCN^{\ot n}\ot\YN[[u\1]].
$$

For each $k=1\lc N$ let $H_k\in\EndCN^{\ot k}$
be the antisymmetrisation map
normalized so that $H_k^{\ts2}=H_k$.
We will make use of the decomposition into an ordered product
in $\EndCN^{\ot k}$
$$
k!\ts(k-1)!\ldots1!\cdot H_k\ts=\prod_{1\leqslant p<k}^\leftarrow
\Bigl(\ts
\prod_{p<q\leqslant k }^\leftarrow\ts
R_{pq}(q-p)
\Bigr)
\Tag{1.4}
$$
where the arrows indicate the order in which the factors are arranged
when the indices $p$ and $q$ increase. Due to this decomposition
the relation \(1.2) implies that
in the algebra $\EndCN^{\ot k}\ot\YN[[u\1]]$ we have
$$
H_k\ot1
\cdot
T_1(u-1)
\ts\ldots\ts
T_k(u-k)
=
T_k(u-k)
\ts\ldots\ts
T_1(u-1)
\cdot
H_k\ot1.
\Tag{1.5}
$$
Here the series in $(u-1)\1,\ts\dots\ts,(u-k)\1$ should be
re-expanded in $u\1$.
\Par
Denote by $\Fk$ the subalgebra in $\EndCN^{\ot k}$ formed by all
elements $X$ such that
$H_k\ts X=H_k\ts X\ts H_k$.
Then by \(1.5) we have
$$
T_1(u-1)
\ts\ldots\ts
T_k(u-k)
\ts\in\ts
\Fk\ot\YN[[u\1]].
$$
We will identify~the~algeb\-ra $\EndLak$ with the subalgebra
in $\EndCN^{\ot k}$ which consists of all the elements
of the form $H_k\hskip1pt X H_k$. Denote by $\phi_k$ the homomorphism
$$
\Fk\to\EndLak:\ts X\mapsto H_k\hskip1pt X.
$$

Let an arbitrary element $Z\in\EndCN$ be fixed. We will now define a
remarkable commutative subalgebra $\ANZ$ in $\YN$. We will call it
$\text{the {\it Bethe subalgebra}}$.
Our definition is only a slight modification of that
given in [KR\ts,\ts Section 2].
\Par
Consider the formal power series $T_1(u)\lc T_N(u)$
with the coefficients in the algebra $\EndCN^{\ot N}\ot\YN$.
The subalgebra $\ANZ$ in $\YN$ is generated by the
coefficients of all the series
$$
B_k(u)
=\tr_N\ot\id
\bigl(\ts
H_N\ot1\cdot
T_1(u-1)
\ts\ldots\ts
T_k(u-k)
\cdot
Z_{k+1}\ldots Z_N\ot1
\ts\bigr)
\Tag{1.8prime}
$$
where $k=1\lc N$.  Let $z_{ij}\in\CC$ be the matrix elements of $Z$;
$$
Z=\sum_{i,j}\ts z_{ij}\ts E_{ij}\ts.
$$
Then in a more conventional notation $B_k(u)$ equals the sum
$$
\sum_{g,\gp}\ts\ts
T_{g(1)\gp(1)}(u-1)
\ldots
T_{g(k)\gp(k)}(u-k)
\cdot
z_{g(k+1)\gp(k+1)}\ldots z_{g(N)\gp(N)}
\cdot\ep\ts/\ts N!
\Tag{1.3}
$$
where $g$ and $\gp$ run through the set of all permutations of $1,2,\dots,N$
while $\ep$ stands for $\sgn\ts g\cdot\sgn\ts\gp$.
Note that the projector $H_N\in\EndCN^{\ot N}$ is one-dimensional.
So by using \(1.5) when $k=N$ we obtain
for the series $B_N(u)$ another expression,
$$
B_N(u)=
\sum_g\ts
T_{g(1),1}(u-1)\ \ldots\ T_{g(N),N}(u-N)
\cdot\operatorname{sgn}g\ts
\Tag{1.99}
$$
where $g$ runs through the set of all permutations of $1,2,\dots,N$.
The series $B_N(u)$ is called the
{\it quantum determinant} for the algebra $\YN$.
The following proposition is well known; its detailed proof can be
found in
[MNO,\ts Section 2].

\proclaim{Proposition 1.1}
The coefficients at $u\1,\ts u^{-2},\ts\dots$ of the series $B_N(u)$
are free generators for the centre of the algebra $\YN$.
\endproclaim

\nt
The proof of the next proposition is also known. Nevertheless,
we will give it here.

\proclaim{Proposition 1.2}
All the coefficients of the series $B_1(u)\lc B_N(u)$ commute.
\endproclaim

\demo{Proof}
By making use of the equality \(1.5) when $k=N$ we obtain that
$$
H_N\ot1
\cdot
T_1(u-1)
\ts\ldots\ts
T_N(u-N)
=
H_N\ot B_N(u).
\Tag{1.6}
$$
The element $T(u)$ belongs to
$$
\id\ot1+\EndCN\ot\YN[[u\1]]\cdot u\1
$$
and is therefore invertible in the algebra $\EndCN\ot\YN[[u\1]]$.
Let $\wT(u)$ be the inverse series. By \(1.6) the series
$B_k(u)$ then equals
$$
B_N(u)\cdot\tr_N\ot\id
\bigl(\ts
H_N\ot1\cdot
\wT_N(u-N)
\ts\ldots\ts
\wT_{k+1}(u-k-1)
\cdot
Z_{k+1}\ldots Z_N\ot1
\ts\bigr).
$$
Put
$$
\wB_k(u)=
\tr_k\ot\id
\bigl(\ts
H_k\ot1\cdot
\wT_k(u-k)
\ts\ldots\ts
\wT_1(u-1)
\cdot
Z_1\ldots Z_k\ot1
\ts\bigr)
\Tag{1.10prime}
$$
where $\wT_1(u)\lc\wT_k(u)$ are regarded as elements of
$\EndCN^{\ot k}\ot\YN[[u\1]]$. Then
$$
B_k(u)=
B_N(u)\ts\wB_{N-k}(u-k)
\cdot
\left(\matrix N\\k\endmatrix\right).
\Tag{1.111111}
$$
By Proposition 1.1 all the coefficients of the series $B_N(u)$ are central
in $\YN$.
Hence it suffices to prove that
$[\wB_k(u),\wB_l(v)]=0$ for all the indices $k,l=1\lc N$.

Let us consider the ordered product
$$
P(u)=\prod_{1\leqslant p\leqslant k}^\rightarrow
\Bigl(\ts
\prod_{k<q\leqslant k+l}^\leftarrow\ts
R_{pq}(u-p+q)
\Bigr)
\in\EndCN^{\ot(k+l)}[u].
\Tag{1.75}
$$
The product $P(u)$ has an inverse in $\EndCN^{\ot(k+l)}(u)$
and commutes with the element $Z^{\ot(k+l)}$.
Moreover,
due to the {\it Yang-Baxter equation} in $\EndCN^{\ot3}[u,v]$
$$
R_{12}(u)\ts R_{13}(u+v)\ts R_{23}(v)=
R_{23}(v)\ts R_{13}(u+v)\ts R_{12}(u)
\Tag{1.8}
$$
and to the decomposition \(1.4) we have
$$
\align
H_k\ot1\cdot P(u)
&=
\prod_{1\leqslant p\leqslant k}^\leftarrow
\Bigl(\ts
\prod_{k<q\leqslant k+l}^\leftarrow\ts
R_{pq}(u-p+q)
\Bigr)\cdot H_k\ot1\ts,
\\
1\ot H_l\cdot P(u)
&=
\prod_{1\leqslant p\leqslant k}^\rightarrow
\Bigl(\ts
\prod_{k<q\leqslant k+l}^\rightarrow\ts
R_{pq}(u-p+q)
\Bigr)\cdot 1\ot H_l\ts.
\endalign
$$
In particular, we have
$$
P(u)\in\Fk\ot\Fl\ts[u].
$$
We will denote
$$
\phi_k\ot\phi_l\bigl(P(u)\bigr)=\Bar P(u).
\Tag{1.755}
$$
The element $\bar P(u)$ has an inverse in
$\EndLak\ot\EndLal(u)$.

By \(1.2) we also have the equality in
$\EndCN^{\ot(k+l)}\ot\YN((u\1,v\1))$
$$
\align
&P(u-v)\ot1
\cdot
T_1(u-1)
\ts\ldots\ts
T_k(u-k)\cdot
T_{k+1}(v-1)
\ts\ldots\ts
T_{k+l}(v-l)
\Tag{1.7}
\\
\qquad\ \qquad
&=T_{k+1}(v-1)
\ts\ldots\ts
T_{k+l}(v-l)
\cdot
T_1(u-1)
\ts\ldots\ts
T_k(u-k)
\cdot
P(u-v)\ot1.
\endalign
$$
Introduce the elements of the algebra
$\EndLak\ot\EndLal\ot\YN[[u\1]]$
$$
\align
K(u)
&=
\phi_k\ot\phi_l\ot\id
\bigl(\ts
\wT_k(u-k)
\ts\ldots\ts
\wT_1(u-1)
\bigr)
\ts,
\\
L(u)
&=
\phi_k\ot\phi_l\ot\id
\bigl(\ts
\wT_{k+l}(u-l)
\ts\ldots\ts
\wT_{k+1}(u-1)
\bigr)
\endalign
$$
and
$
W=
{\Lambda\hskip-2pt}^{\raise0.5pt\hbox{$\scriptstyle k$}}\hskip-1ptZ
\ot
{\Lambda\hskip-2pt}^{\raise0.5pt\hbox{$\scriptstyle l$}}\hskip-1ptZ
\ot1.
$
Then we have
$$
\wB_k(u)\ts\wB_l(v)
=
\tr\ot\id
\bigl(\ts
K(u)\ts L(v)\cdot W
\bigr)
$$
where $\tr$ stands for the restriction of $\tr_{k+l}$ onto
$\EndLak\ot\EndLal$.
But by applying the homomorphism $\phi_k\ot\phi_l\ot\id$ to \(1.7) we get
the equalities
$$
\align
K(u)\ts L(v)\cdot W\cdot\bar P(u-v)\ot1
&=
K(u)\ts L(v)\cdot\bar P(u-v)\ot1\cdot W
\\
&=
\bar P(u-v)\ot1
\cdot L(v)\ts K(u)\cdot W.
\endalign
$$
Therefore we obtain that
$$
\tr\ot\id
\bigl(\ts
K(u)\ts L(v)\cdot W
\bigr)
=
\tr\ot\id
\bigl(\ts
L(v)\ts K(u)\cdot W
\bigr)
=
\wB_l(v)\ts\wB_k(u)\quad\square
$$
\enddemo

\proclaim{Theorem 1.3}
Suppose that the element $Z\in\EndCN$ has a simple spectrum.
Then the subalgebra $\ANZ$ in $\YN$ is maximal commutative.
The coefficients at $u\1,\ts u^{-2},\ts\dots$ of the series
$B_1(u)\lc B_N(u)$
are free generators for $\ANZ$.
\endproclaim

\nt
The proof of this theorem will be given in Section 2.
We will end up this section with comparing our definition
of the Bethe subalgebra in $\YN$ with that
given in [KR\ts,\ts Section 2].
The defining relations \(1.2) show that the assignment
$T(u)\mapsto\wT(u)$ defines an antiautomorphism of the
algebra $\YN$. Denote by $\eta$ this antiautomorphism.
It then follows from \(1.6) that
$$
\eta\bigl(B_N(u)\bigr)=B_N(u)\1.
$$
Furthermore, we have
$$
\eta^2\bigl(T(u)\bigr)=T(u+N)\cdot B_N(u)/B_N(u+1);
\nopagebreak
$$
the proof of the latter statement can be found
for instance in [MNO\ts,\ts Section 5].

The commutative subalgebra of $\YN$ considered in [KR\ts,\ts Section 2] is
generated by the coefficients of all the series
$$
\eta\1\bigl(\wB_k(u)\bigr)
=
\tr_k\ot\id
\bigl(\ts
H_k\ot1\cdot
T_1(u-1)
\ts\ldots\ts
T_k(u-k)
\cdot
Z_1\ldots Z_k\ot1
\ts\bigr)
\Tag{1.7777777}
$$
where $k=1\lc N$.
This commutative subalgebra is not maximal if the element
$Z\in\EndCN$ is not invertible. In that case we have $\wB_N(u)=0$.
However, due to \(1.111111)
the commutative subalgebra in $\YN$ generated by the coefficients of
all the series \(1.7777777) along with $B_N(u)$ coincides with
$\eta\1\bigl(\ANZ\bigr)$. The latter subalgebra is maximal
for any element $Z$ with a simple spectrum by Theorem 1.3.

\Section{Proof of Theorem 1.3}

\nt
We will reduce the proof to several lemmas.
Some of them are rather general and will be used again in Section 4.
We will employ methods from [C3] and [K1].

Let us equip the algebra $\YN$ with an ascending filtration
by setting degrees
\line{of its generators as
$\deg T_{ij}^{(r)}=r$. The linear subspace in $\YN$ consisting of all
the} elements with degrees not greater than $r$ will be denoted by
$\YrN$. Consider the graded algebra
$$
\YgN=\underset{r\geqslant0}\to\oplus\ts\YrN/\YrmN
$$
corresponding to the filtered algebra $\YN$.
The defining relations \(1.11) show that the algebra
$\YgN$ is commutative.

There is a natural Poisson algebra structure on $\YgN$.
For any two elements
$X,Y$ in $\YN$ of degrees $p,q$ respectively the Poisson
bracket of their images $x,y$ in $\YgN$ is
$$
\{\ts x,y\ts\}=[\ts X,Y\ts]\ \operatorname{mod}\ts\YpqN.
$$
Let $\AgNZ$ be the image in $\YgN$ of the Bethe subalgebra
$\ANZ$ in the Yangian $\YN$.
To prove the first statement of Theorem 1.3
it suffices to show that the
subalgebra $\AgNZ$ in $\YgN$ is maximal involutive.
\Par
We will denote by $t_{ij}^{(r)}$ the image in $\YgN$
of the generator $T_{ij}^{(r)}$ of the \text{algebra} $\YN$.
All the elements $t_{ij}^{(r)}\hskip-1pt$ are free generators of the
commutative \text{algebra}~$\hskip-1pt\YgN\hskip-1pt$.
The proof of this assertion
can be found for instance in [MNO, Section 1].
We will identify $\YgN$ with the
symmetric algebra $\operatorname{S}(\glN[\hskip1pt t \hskip1pt])$
over the {polynomial current Lie algebra} $\glN[\hskip1pt t \hskip1pt]$.
The generator $t_{ij}^{(r)}$
will be identified with the element $E_{ij}\ts t^{\ts r-1}$.
We will also set $t_{ij}^{(0)}=\de_{ij}\cdot1$.
Then by (1.11)
$$
\{\ts t_{ij}^{(p)},t_{kl}^{(q)}\ts\}=\ts
\sum_{r=1}^{\min(p,q)}
\bigl(\ts
t_{kj}^{(r-1)}\ts t_{il}^{(p+q-r)}-
t_{kj}^{(p+q-r)}\ts t_{il}^{(r-1)}
\ts\bigr);
\ \quad
p,q=1,2,\ldots\ts\ts.
\Tag{2.1}
$$

Consider the formal power series in $u\1$
with the coefficients in $\YgN$
$$
t_{ij}(u)=
t_{ij}^{(0)}+t_{ij}^{(1)}\ts u\1+t_{ij}^{(2)}\ts u^{-2}+\ldots\ts.
$$
Then the image in the algebra $\YgN[[u\1]]$ of the series
$T_{ij}(u-h)\in\YN[[u\1]]$ equals $t_{ij}(u)$
for any $h\in\CC$. Therefore the image of the series
$B_k(u)$ in the algebra $\YgN[[u\1]]$ equals
$$
\sum_{g,\gp}\ts\ts
t_{g(1)\gp(1)}(u)
\ts\ldots\ts
t_{g(k)\gp(k)}(u)
\cdot
z_{\ts g(k+1)\gp(k+1)}\ldots\ts z_{\ts g(N)\gp(N)}
\cdot\ep\ts/\ts N!
\Tag{2.11111111}
$$
where $g,\gp$ and $\ep$ are the same as in \(1.3).
Note that the
coefficients of the latter series are homogeneous in $\YgN$.
\Par
For each $M=1,2\ldots$ let us denote by $\YgnN$ the ideal in the algebra
$\YgN$ generated by all the elements $t_{ij}^{(r)}$ where $r\geqslant M$.
Note that by \(2.1) we have
$$
\{\ts t_{ij}^{(p)},t_{kl}^{(q)}\ts\}\ts\in\ts\YgnN,
\ \quad
M=\max(p,q).
$$
Therefore each $\YgnN$ is a Poisson ideal in $\YgN$. Due to the
following general lemma it suffices to prove that for each $M$
the image of $\AgNZ$ in the
quotient Poisson algebra $\YgN/\YgnN$ is maximal involutive.

Let $\operatorname{X}$ be any graded Poisson algebra and
$
\operatorname{X}_1\supset \operatorname{X}_2\supset\ldots
$
be a descending chain of Poisson ideals in $\operatorname{X}$ such that
for each $M=1,2,\ldots$ the ideal $\operatorname{X}_M$ contains only
elements of degrees not less than $M$. Let $b_1,b_2,\ldots$ be a sequence
of homogeneous elements of $\operatorname{X}$ in involution.

\proclaim{Lemma 2.1}
Supose that for each $M=1,2,\ldots$
the images of $b_1,b_2,\ldots$ in the Poisson algebra
$\operatorname{X}/\operatorname{X}_M$ generate
a maximal involutive subalgebra. Then $b_1,b_2,\ldots$ generate
a maximal involutive
subalgebra in $\operatorname{X}$.
\endproclaim

\demo{Proof}
Fix a homogeneous element $x\in\operatorname{X}$
in involution with each of $b_1,b_2,\ldots$ and set $M=1+\deg x$.
By our assumption
$
x=f(b_1,b_2,\ldots)+y
$
for certain polynomial $f$ and some
$y\in\operatorname{X}_M$.
Then $x$ is the sum of the terms in
$f(b_1,b_2,\ldots)$ of degrees smaller than $M\ \ \square$
\enddemo

\nt
By definition, we have $\YgoN=\operatorname{S}(\glN)$. From now we shall keep
$M\geqslant1$ fixed. For each $r=1\lc M$
we will denote by $x_{ij}^{(r)}$
the image of the generator $t_{ij}^{(r)}$ in the quotient
Poisson algebra $\YgN/\YgnoN$. By \(2.1) in the latter algebra we have
$$
\{\ts x_{ij}^{(p)},x_{kl}^{(q)}\ts\}=\ts
\sum_{r=\max(1,\ts p+q-M)}^{\min(p,q)}
\bigl(\ts
x_{kj}^{(r-1)}\ts x_{il}^{(p+q-r)}-
x_{kj}^{(p+q-r)}\ts x_{il}^{(r-1)}
\ts\bigr)
\Tag{2.11}
$$
where $x_{ij}^{(0)}=\de_{ij}\cdot1$. Consider the vector subspace
in the
Lie algebra $\glN[\hskip1pt t \hskip1pt]$
$$
\uMN=\glN
+\ldots+\glN\cdot t^{\ts M-1}.
\Tag{2.2222}
$$
As a commutative algebra
the quotient $\YgN/\YgnoN$ can be identified with the symmetric algebra
$\operatorname{S}(\uMN\!)$. Further, we will identify
$\YgN/\YgnoN$
with the algebra $\operatorname{P}(\uMN\!)$ of the polynomial functions
on $\uMN$.
The generator $x_{ij}^{(r)}$ will be regarded as the coordinate
function corresponding to the vector $E_{ij}\cdot t^{\ts r-1}$.
\Par
Let us now make use of the following observation.
The relations \(1.2) imply that for any invertible element
$G\in\EndCN$ the assignment
$$
T(u)\mapsto G\ot1\cdot T(u)\cdot G\1\ot1
$$
determines an automorphism of the algebra $\YN$.
The image of the subalgebra
$\ANZ$ with respect to this automorphism coincides with
$\operatorname{B}(\glN, G\1\hskip1pt Z\hskip1pt G)$. Hence
it suffices to assume that $z_{ij}=z_i\cdot\de_{ij}$ where all
$z_i\in\CC$ are pairwise distinct.
{}From now on until the end of this section we shall keep to this assumption.

Let $\bN$ be the Borel subalgebra in $\glN$ spanned by the elements
$E_{ij}$ with $i\leqslant\nomathbreak j$. Fix any
principal nilpotent element of the opposite Borel subalgebra of the form
$$
E=\ep_1\ts E_{21}+\ldots+\ep_{N-1}\ts E_{N,N-1}
$$
where $\ep_1\lc\ep_{N-1}\neq0$.
We will consider the Poisson bracket \(2.11) on
$\operatorname{P}(\uMN\!)$
in a neighbourhood of the point
$$
E^{(M)}=E\cdot t^{\ts M-1}\in\uMN.
$$
Introduce the affine subspace in $\uMN$
$$
\kN=E^{(M)}+\uMN\cap\bN[\hskip1pt t \hskip1pt].
$$

For each $k=1\lc N$ denote by $b_k(u)$ the image of
$B_k(u)\in\YN[[u\1]]$ in
$$
\YgN/\YgnoN[[u\1]]=
\operatorname{P}(\uMN\!)[[u\1]].
$$
This image is in fact a polynomial in $u\1$
of the degree $k\hskip1pt M$. We will write
$$
b_k(u)=b_k^{(0)}+b_k^{(1)}\ts u\1+\ts\ldots\ts+b_k^{(kM)}\ts u^{-kM}\ts.
$$
Here all
the coefficients $b_k^{(r)}$ are in involution due to Proposition 1.2.
Consider them as polynomial functions on the vector space $\uMN$.
Of the proof of Theorem 1.3 the next lemma is the main part;
cf. [K1\ts,\ts Section 4] and [RT\ts,\ts Section 4]. Another approach to the
proof of this theorem was described in [RS\ts,\ts Section 3].

\proclaim{Lemma 2.2}
The restrictions of the functions $b_k^{(r)}$ onto the affine
subspace $\kN$
generate the whole algebra of polynomial functions on $\kN$.
\endproclaim

\demo{Proof}
Introduce the polynomial in $u$ of the degree $M-1$
$$
x_{ij}(u)=
x_{ij}^{(1)}\ts u^{M-1}+x_{ij}^{(2)}\ts u^{M-2}+\ts\ldots\ts+
x_{ij}^{(M)}.
$$
Then
$$
\align
u^{kM}\ts b_k(u)=
\sum_{g,h}\ts\ts
\bigl(u^{M}\delta_{g(1)h(1)}+x_{g(1)h(1)}(u)\bigr)
&\ldots
\bigl(u^{M}\delta_{g(k)h(k)}+x_{g(k)h(k)}(u)\bigr)
\\
\hskip5.3cm
\times\ts\ts
z_{\ts g(k+1)h(k+1)}
&\ldots\ts
z_{\ts g(N)h(N)}
\cdot\varepsilon\ts/\ts N!
\endalign
$$
where $g,h$ and $\varepsilon$ are the same as in \(1.3) and \(2.11111111).
Denote by $\hskip-1pt X(u)\hskip-1pt$
the square matrix of order $\hskip-1pt N$ formed by all the
polynomials $x_{ij}(u)$.
Then we have the Laplace expansion
$$
\det\bigl(u^M+X(u)+Z\ts v\bigr)=
v^N\det Z
+\sum_{k=1}^N\ts\ts
u^{kM}\ts b_k(u)\ts v^{\,N-k}
\left(\!\matrix N\\k\endmatrix\!\right)
\Tag{2.*}
$$
where $Z$ is now regarded as a diagonal matrix of order $N$ with
pairwise distinct diagonal entries $z_1\lc z_N\in\CC$. Denote by $F(u,v)$
the polynomial in $u,v$ obtained by restricting the coefficients of
\(2.*) onto $\kN$.
\Par
On the subspace $\kN$ we have $x_{i+1,i}(u)=\ep_i$
and $x_{ij}(u)=0$ if $i-j>1$.
Take the functions $x_{i,i+d}^{(r)}$ with $d=0\lc N-1$
and $r=1\lc M$ as coordinates
on $\kN$. Endow the set of the pairs $(d,r)$ with the lexicographical order.
We will prove consecutively for each of these pairs that
$x_{i,i+d}^{(r)}$ are polynomials in the coefficients of $F(u,v)$.
\Par
Assume that $\deg v=M$ while $\deg u=1$. Consider the terms of $F(u,v)$
with the total degree $M(N-d)-r$. Their sum has the form
$$
f(u,v)+(-1)^{d}\cdot
\sum_{i=1}^{N-d}\ts\ts
x_{i,i+d}^{(r)}\ts\ts u^{M-r}
\cdot g_i(u^M,v)\ts\ts\ep_i\hskip1pt\ldots\hskip1pt\ep_{i+d-1}
$$
where
$$
g_i(u,v)=
\prod_{j\neq i\lc i+d}\ts (u+z_j\ts v).
$$
Here the coefficients of the polynomial $f(u,v)$
depend only on $z_1\lc z_N$ and
$x_{ij}^{(s)}$ with the pair $(j-i,s)$ preceding $(d,r)$.
It now remains to observe that all the $N-d$
polynomials
$g_i(u,v)$ are linearly independent. Indeed, 
for each $i=1\lc N-d$ we have
$$
g_i(-z_i\ts,1)=\prod_{j\neq i\lc i+d}\ts (z_j-z_i)\neq0
$$
while for $1\leqslant i<j\leqslant N-d$ we have
the equality $g_j(-z_i\ts,1)=0$
\enddemos

\nt
Note that here $b_k^{(0)}\in\CC$ for any $k=1\lc N$.
The total number of the remaining coefficients $b_k^{(r)}$
with $r=1\lc k\hskip1pt M$ equals
$$
M+2M+\ldots+N\hskip1pt M=M\hskip1pt N\hskip1pt(N+1)/2=\dim\ts\kN\ts.
$$

\proclaim{Corollary 2.3}
All the coefficients $b_k^{(r)}$ with $k=1\lc N$ and
$r=1\lc k\hskip1pt M$ are algebraically independent.
\endproclaim

\nt
This corollary is valid for any $M\geqslant1$ and therefore already
implies the second statement of Theorem 1.3. Namely, all
the coefficients at $u\1,\ts u^{-2},\ts\dots$ of the series
$B_1(u)\lc B_N(u)$ are free generators for $\ANZ$. Denote
$$
D=\dim\ts\uMN-\dim\ts\kN=M\hskip1pt N\hskip1pt(N-1)/2\ts.
$$

\proclaim{Lemma 2.4}
There is an open neighbourhood of the point
$E^{(M)}$ in $\uMN$
where the rank of the Poisson bracket \(2.11) equals $2D$.
\endproclaim

\demo{Proof}
The elements $x_{ij}^{(1)}$
generate
a Poisson subalgebra isomorphic to the algebra
$\operatorname{P}(\glN)$ of the polynomial functions on $\glN$
with the standard Poisson bracket:
$$
\{\ts x_{ij}^{(1)},x_{kl}^{(1)}\ts\}=\ts
\de_{kj}\ts x_{il}^{(1)}-
\de_{il}\ts x_{kj}^{(1)}.
\Tag{2.4444}
$$
By \(2.11) for any $p,q\geqslant1$ the value of the function
$\{\ts x_{ij}^{(p)},x_{kl}^{(q)}\ts\}$
at the point $E^{(M)}$ equals
$$
\de_{p+q,M+1}
\cdot
(\ts\de_{kj}\ts\de_{i,l+1}\ts\ep_l-\de_{il}\ts\de_{k,j+1}\ts\ep_j\ts).
$$
Therefore the rank of the bracket \(2.11) at this point is
$M$ times the rank of the standard Poisson bracket on
$\operatorname{P}(\glN)$ in the point $E$. The latter rank is
$N^2-N$. So there exists
an open neighbourhood of the point $E^{(M)}$ in
$\uMN$ where
the rank of the Poisson bracket \(2.11) is not smaller than
$
M\hskip1pt (N^2-N)=2D.
$
Due to Proposition 1.1 and Corollary 2.3
there are $MN$ algebraically independent elements
$
b_N^{(1)}\lc b_N^{(MN)}
$
in the
centre of the Poisson algebra
$\operatorname{P}(\uMN\!)$.
So the rank of the Poisson bracket cannot exceed
$M\hskip1pt(N^2-N)=2D$ at any point $\ \square$
\enddemo

\nt
Let us now fix any
$x\in\operatorname{P}(\uMN\!)$
in involution with all the elements $b_k^{(r)}$. To complete the
proof of Theorem 1.3 we have to demonstrate that the element $x$
is then a polynomial in $b_k^{(r)}$. For any collection
$$
H=\bigl(\ts h_k^{(r)}\in\CC\ |\ k=1\lc N;\ r=1\lc kM\ts\bigr)
$$
denote by $\Cal S_H$ the subset in $\uMN$
where the values of the functions $b_k^{(r)}$ are $h_k^{(r)}$
respectively.
Due to Lemma 2.2 there is an open neighbourhood $\Cal E$ of the point
$E^{(M)}\in\uMN$ such that
every non-empty intersection $\Cal E\cap\Cal S_H$ is
transversal to $\kN$.
Again due to Lemma 2.2 to demonstrate that $x$
is a polynomial in $b_k^{(r)}$ it now
suffices to prove the following statement; cf. [C3].

\proclaim{Lemma 2.5}
One can choose the open neighbourhood $\Cal E$
so that the function
$x$ is constant on every intersection $\Cal E\cap\Cal S_H$.
\endproclaim

\demo{Proof}
For any polynomial $b$ in $b_k^{(r)}$ consider
the respective flow in $\uMN$
$$
{d\hskip0.5pt x_{ij}^{(r)}}
/
{d\hskip1pt t}
=
\{\ts b,x_{ij}^{(r)}\ts\}\ts;
\qquad
i,j=1\lc N;
\ \
r=1\lc M
\Tag{2.2}
$$
where $t$ stands for the coordinate on the line $\RR$.
Since $\{\ts b,x\ts\}=0$
the function $x$ is constant along every trajectory of this flow.
For any point $F\in\uMN$ denote by $\Cal T_F$
the collection of the trajectories of the flows \(2.2) passing
through $F$ for all polynomials $b$ in $b_k^{(r)}$.
Due to Lemma~2.2 and to Lemma 2.4 by the Liouville theorem we can choose
the open neighbourhood $\Cal E$ of $E^{(M)}$
so that
$$
F\in S_H
\quad\Rightarrow\quad
\Cal T_F\cap\Cal E=\Cal S_H\cap\Cal E
\qquad\square
$$
\enddemo

Now we will consider
the commutative subalgebra $\pi\bigl(\ANZ\bigr)$ in $\UN$.
There is a canonical filtration on the algebra $\UN$.
Consider the involutive subalgebra in the graded algebra  
$\operatorname{S}(\glN)$ corresponding to $\pi\bigl(\ANZ\bigr)$. 
We will identify the Poisson algebras
$\operatorname{S}(\glN)$ and $\operatorname{P}(\glN)$;
the element $E_{ij}\in\glN$ will be identified with the respective
coordinate function.
Let $X$ be
the square matrix of order $N$ formed by these
functions. The subalgebra  in $\operatorname{P}(\glN)$
corresponding to $\pi\bigl(\ANZ\bigr)$
is generated by the coefficients of the polynomial in $u\ts,v$
$$
\det\bigl(u+X+Z\ts v\bigr).
$$
It is well known that the coefficients of this polynomial
at $u^k\ts v^l$ with $k+l<N$ are algebraically independent
for any $Z$ with a simple spectrum;
see [H\ts,\ts Section 2] and [MF\ts,\ts Section 4].
Lemma 2.2 at $M=1$ provides an elementary proof of this fact.
Moreover, it shows that these coefficients then
generate a maximal involutive subalgebra in $\operatorname{P}(\glN)$.
Thus we obtain the following statement.

\proclaim{Proposition 2.6}
\!The subalgebra $\pi\bigl(\ANZ\bigr)$\hskip-1pt in \hskip-1pt$\UN$
\hskip-1pt is \text{maximal commutative.}
\endproclaim

\nt
The analogues of this proposition for
the universal enveloping algebras
$\operatorname{U}(\hskip1pt\so_N)$ and $\operatorname{U}(\hskip1pt\sp_N)$
of the other classical Lie algebras will be given
in Section 4.

\Section{Bethe subalgebras in twisted Yangians}

\nt
In the previous section we considered only the Yangian of the
Lie algebra $\glN$. We~will start this section with describing
analogues of this Yangian for the other classical Lie algebras
$\soN$ and $\spN$;\ts\ in the latter case $N$ has to be even. These
analogues have been introduced in [O2]; see also [MNO\ts,\ts Section 3].
Then we will construct the
respective analogues of the Bethe subalgebra $\ANZ$~in~$\YN$.

Let $\aN$ be one of the classical Lie algebras $\soN$ and $\spN$.
We will regard $\aN$ as an involutive subalgebra in $\glN$.
Let $\sigma$ be the corresponding involutive automorphism of the
Lie algebra $\glN$.
The superscript $^\prime$ will denote transposition in $\EndCN$ with
respect to the symmetric or alternating bilinear form on $\CC^N$
preserved by the subalgebra $\aN$ in $\glN$. As well as above
$E_{ij}\in\EndCN$ will be the standard matrix units.
But now we will let the indices $i$ and $j$ run through the set
$\{\ts-n\lc-1,1\lc n\ts\}$ if $N=2n$ and the set
$\{\ts-n\lc-1,0,1\lc n\ts\}$ if $N=2n+1$.

Put $\ep_{ij}=\sgn\ts i\cdot \sgn\ts j$ if $\aN=\spN$ and
$\ep_{ij}=1$ if $\aN=\soN$.
We will choose the symmetric or alternating bilinear form on $\CC^N$
so that
$$
E_{ij}^{\,\prime}=\ep_{ij}\cdot E_{-j,-i}\ts.
$$
If we regard $E_{ij}$ as generators of the universal enveloping
algebra $\UN$ then $\sigma(E_{ij})=-E_{ij}^{\,\prime}$ so that
$F_{ij}=E_{ij}-E_{ij}^{\,\prime}$ are generators of the
algebra $\UaN$.

Let us now introduce the element of the algebra $\EndCN\ot\YN[[u\1]]$
$$
\tT(u)=\sum_{i,j}\ts E_{ij}^{\,\prime}\ot T_{ij}(u)
$$
and consider the series with the coefficients in $\EndCN\ot\YN$
$$
\gather
S(u)=T(u)\ts\tT(-u)=\sum_{i,j}\ts E_{ij}\ot S_{ij}(u),
\Tag{3.01}
\\
\tS(u)=\sum_{i,j}\ts E_{ij}^{\,\prime}\ot S_{ij}(u).
\endgather
$$
Then
$$
S_{ij}(u)=
\de_{ij}\cdot1+S_{ij}^{(1)}\ts u\1+S_{ij}^{(2)}\ts u^{-2}+\ldots
$$
for certain elements $S_{ij}^{(1)},S_{ij}^{(2)},\ldots\in\YN$.
By definition, the {\it twisted Yangian} $\YaN$ is the subalgebra in $\YN$
generated by all the elements $S_{ij}^{(r)}$.
This definition along with \(1.*) implies that $\YaN$ is
a left coideal in the Hopf algebra $\YN$:
$$
\Delta\bigl(\YaN\bigr)\subset\YN\ot\YaN.
$$
Introduce also the element of $\EndCN^{\ot2}[u]$
$$
\tR(u)=
u\cdot\id-\sum_{i,j}\ts E_{ij}^{\,\prime}\ot E_{ji}=
u\cdot\id-\sum_{i,j}\ts E_{ij}\ot E_{ji}^{\,\prime}\ts.
$$
Later on we will employ the eqiality
$$
\tR(u)\ts\tR(-u+N)=(Nu-u^2)\cdot\id.
\Tag{3.0}
$$
The relation \(1.2) implies that in the algebra
$\EndCN^{\ot2}\!\ot\YN\ts((u\1,v\1))$
$$
\tT_1(u)\cdot\tR(u-v)\ot1\cdot T_2(v)=
T_2(v)\cdot\tR(u-v)\ot1\cdot\tT_1(u).
\Tag{3.1}
$$
By making use of the latter relation and by applying \(1.2) again we
obtain the relation in
$\EndCN^{\ot2}\!\ot\YaN\ts((u\1,v\1))$
$$
\align
&R(u,v)\ot1\cdot S_1(u)\cdot\tR(-u-v)\ot1\cdot S_2(v)
\Tag{3.2}
\\
&=S_2(v)\cdot\tR(-u-v)\ot1\cdot S_1(u)\cdot R(u,v)\ot1.
\endalign
$$

\nt
The proof of the following proposition is contained in [MNO\ts,\ts Section 3].
Whenever in the present article the double sign $\pm$ or $\mp$ occurs,
the upper one will correspond to the case $\aN=\soN$ while the lower sign
will correspond to the case $\aN=\spN$.

\proclaim{Proposition 3.1}
In the algebra $\EndCN\ot\YaN\ts[[u\1]]$ we have the relation
$$
S(u)-\tS(-u)=\mp\ts\bigl(S(u)-S(-u)\bigr)/\ts2u\ts.
\Tag{3.3}
$$
The relations \(3.2) and \(3.3) yield the defining relations
for the generators
$S_{ij}^{(r)}$ of the subalgebra $\YaN$ in $\YN$.
\endproclaim

\nt
In a more conventional notation \(3.2) and \(3.3)
can be rewritten respectively as the collections of relations
for all possible indices $i,j,k,l$
$$
\align
(u^2-v^2)\cdot[S_{ij}(u),S_{kl}(v)]=
(u+v)\cdot
\bigl(
S_{kj}(u)\ts S_{il}(v)&-S_{kj}(v)\ts S_{il}(u)
\bigr)
\\
-(u-v)\cdot
\bigl(
\ep_{k,-j}\,S_{i,-k}(u)\,S_{-j,l}(v)&-\ep_{i,-l}\,S_{k,-i}(v)\,S_{-l,j}(u)
\bigr)
\\
+\ts\ep_{i,-j}\cdot
\bigl(
S_{k,-i}(u)\,S_{-j,l}(v)&-S_{k,-i}(v)\,S_{-j,l}(u)
\bigr)
\Tag{3.22}
\endalign
$$
and for all $i,j$
$$
\hskip0.5in
S_{ij}(u)-\ep_{ij}\ts S_{-j,-i}(-u)=\mp
\bigl(S_{ij}(u)-S_{ij}(-u)\bigr)/\ts2u\ts.
\Tag{3.2121}
$$
The algebra $\hskip-1pt\YaN$ contains the universal enveloping algebra
$\hskip-1pt\UaN\hskip-1pt$
as~a~sub\-algebra; due to \(3.22)
the assignment $F_{ij}\mapsto S_{ij}^{(1)}$
defines the embedding. Moreover, due to \(3.22) and \(3.2121)
we obtain the following corollary to Proposition 3.1.

\proclaim{Corollary 3.2}
There is a homomorphism
$$
\rho:\ts\YaN\to\UaN:\ts\ts S_{ij}(u)\mapsto \de_{ij}+F_{ij}\ts (u\pm1/2)^{-1}.
\Tag{3.1005}
$$
\endproclaim

\nt
The homomorphism $\rho$ is by definition identical on the subalgebra
$\UaN$. We will regard \(3.1005) as an analogue of the evaluation
homomorphism \(1.1005).

The Yang-Baxter equation \(1.8) in $\EndCN^{\ot3}[u,v]$
along with \(1.0) implies
$$
\align
R_{12}(u)\ts\tR_{13}(v)\ts\tR_{23}(u+v)
&=
\tR_{23}(u+v)\ts\tR_{13}(v)\ts R_{12}(u),
\Tag{3.81}
\\
R_{13}(u)\ts\tR_{12}(v)\ts\tR_{23}(u+v)
&=
\tR_{23}(u+v)\ts\tR_{12}(v)\ts R_{13}(u),
\Tag{3.82}
\\
R_{23}(u)\ts\tR_{12}(v)\ts\tR_{13}(u+v)
&=
\tR_{13}(u+v)\ts\tR_{12}(v)\ts R_{23}(u).
\Tag{3.83}
\endalign
$$
For each $k=1\lc N$ consider the element of
$\EndCN^{\ot k}\ot\YaN[[u\1]]$
$$
S(u,k)=
\omega(u)\cdot
\prod_{1\leqslant p\leqslant k}^\rightarrow
\Bigl(
S_p(u-p)
\cdot
\prod_{p<q\leqslant k}^\rightarrow
\tR_{pq}(p+q-2u)\ot1
\Bigr)
$$
where
$$
\omega(u)=
\prod_{1\leqslant p<q\leqslant k}
(p+q-2u)\1\ts.
$$
Due to the decomposition \(1.4)
by applying the relation \(3.2) repeatedly along with
\(3.81) to \(3.83) we obtain the equality
in the algebra $\EndCN^{\ot k}\ot\YaN((u\1))$
$$
H_k\ot1\cdot S(u,k)=
\omega(u)
\cdot\!
\prod_{1\leqslant p\leqslant k}^\leftarrow
\Bigl(
\prod_{p<q\leqslant k}^\leftarrow
\tR_{pq}(p+q-2u)\ot1
\cdot S_p(u-p)
\Bigr)
\cdot H_k\ot1
\ts;
$$
it is an analogue of the equality \(1.5) in $\EndCN^{\ot k}\ot\YN[[u\1]]$.
In particular,
$$
S(u,k)\ts\in\ts\Fk\ot\YaN[[u\1]].
\Tag{3.5}
$$

Let us keep fixed an arbitrary element $Z\in\EndCN$
but now assume that $Z^{\,\prime}=Z$ or $Z^{\,\prime}=-Z$.
We will now construct a certain commutative subalgebra in $\YaN$.
It will be regarded as an analogue of the Bethe
subalgebra in the Yangian $\YN$ and denoted by $\AaNZ$.

For each $k=1\lc N$ introduce the element of $\EndCN^{\ot k}[[u\1]]$
$$
Z(u,k)=
\omega(u)
\cdot
\prod_{1\leqslant p\leqslant k}^\rightarrow
\Bigl(
Z_p\ts
\cdot
\prod_{p<q\leqslant k}^\rightarrow
\tR_{pq}(p+q-2u)
\Bigr)\ts.
\Tag{3.007095}
$$
Observe that due to the assumption
that  $Z^{\,\prime}=Z$ or $Z^{\,\prime}=-Z$ a relation
similar to \(3.2) holds in the algebra $\EndCN^{\ot2}[u,v]$.
It can be verified directly that
$$
R(u,v)\ts Z_1\ts\tR(-u-v)\ts Z_2=
Z_2\ts\tR(-u-v)\ts Z_1\ts R(u,v).
\Tag{3.25}
$$
By applying the relation \(3.25) repeatedly along with the equalities
\(3.81) to \(3.83) we obtain that
$$
H_k\cdot Z(u,k)=
\omega(u)
\cdot\!
\prod_{1\leqslant p\leqslant k}^\leftarrow
\Bigl(
\prod_{p<q\leqslant k}^\leftarrow
\tR_{pq}(p+q-2u)
\cdot Z_p
\Bigr)
\cdot H_k
\ts.
$$
In particular, we obtain that $Z(u,k)\ts\in\ts\Fk(u)$.
Introduce
the formal power series in $u\1$ with the coefficients
in $\YaN$
$$
A_k(u)=
\tr_N\ot\id
\bigl(
H_N\ot1
\cdot
S_{1\lc k}(u,k)
\cdot
I(u)\ot1
\cdot
Z_{k+1\lc N}(u+N/2-k,N-k)\ot1
\bigr)
$$
where
$$
I(u)=
\prod_{1\leqslant p\leqslant k}^\rightarrow
\Bigl(
\prod_{k<q\leqslant N}^\rightarrow
\tR_{pq}(p+q-2u)/(p+q-2u)
\Bigr).
$$

By definition,
the subalgebra $\AaNZ$ in the twisted Yangian $\YaN$ is generated
by all coefficients of the series
$A_1(u)\lc A_N(u)$.
By using \(3.5) when $k=N$ we get the analogue of \(1.6)
$$
H_N\ot1\cdot S(u,N)=H_N\ot A_N(u).
\Tag{3.6}
$$
The proof of the next proposition is also contained in
[MNO\ts,\ts Section 3].
Let us put $\theta(u)=1+N/(1-2u)$ if $\aN=\spN$. If $\aN=\soN$
we put $\theta(u)=1$.

\proclaim{Proposition 3.3}
We have the equality in the algebra $\YN[[u\1]]$
$$
A_N(u)\ts\theta(u)=B_N(u)\ts B_N(N-u+1).
$$
The coefficients at $u^{-2},\ts u^{-4},\ts\dots$ of the series
$A_N(u)$
are free generators for the centre of the algebra $\YaN$.
\endproclaim

\nt
The series $A_N(u)\ts\theta(u)$ is called the {\it quantum determinant}
for the algebra $\YaN$,
or {\it Sklyanin determinant} [MNO].
There is an analogue of the expression \(1.99) for this series,
it has been proposed in [M].
Now we will prove the following theorem.

\proclaim{Theorem 3.4}
All the coefficients of the series $A_1(u)\lc A_N(u)$ commute.
\endproclaim

\demo{Proof}
We will employ arguments already used in the proof of Proposition 1.3.
We shall keep to the notation introduced therein.
We will also use a method from [S\ts,\ts Section 3].
The element $S(u,k)$ belongs to
$$
\id\ot1+\EndCN^{\ot k}\ot\YaN[[u\1]]\cdot u\1
$$
and is therefore invertible in the algebra $\EndCN^{\ot k}\ot\YaN[[u\1]]$.
Let $\wS(u,k)$ be the inverse series. Put
$$
\wA_k(u)=
\tr_k\ot\id
\bigl(\ts
H_k\ot1
\cdot
\wS(u,k)
\cdot
Z(u+N/2\ts,k)\ot1
\ts\bigr).
\Tag{3.333}
$$
Then by \(3.6)
$$
A_k(u)=
A_N(u)\ts\wA_{N-k}(u-k)\cdot
\left(\matrix N\\k\endmatrix\right).
\nopagebreak
$$
By Proposition 3.3 all the coefficients of the series $A_N(u)$ are central
in $\YaN$.
Hence it suffices to prove that
$[\wA_k(u),\wA_l(v)]=0$ for all the indices $k,l=1\lc N$.

Let $P(u)$ denote the same ordered product as in \(1.75).
Consider also the product
$$
Q(u)=
\prod_{1\leqslant p\leqslant k}^\rightarrow
\Bigl(\ts
\prod_{k<q\leqslant k+l}^\rightarrow\ts
\tR_{pq}(u+p+q)
\Bigr)
\in\EndCN^{\ot(k+l)}[u].
\Tag{3.75}
$$
The product $Q(u)$ has an inverse in $\EndCN^{\ot(k+l)}(u)$.
Due to the decomposition \(1.4) and the equations \(3.81) to \(3.83)
we have
$$
\align
H_k\ot1\cdot Q(u)
&=
\prod_{1\leqslant p\leqslant k}^\leftarrow
\Bigl(\ts
\prod_{k<q\leqslant k+l}^\rightarrow\ts
R_{pq}(u+p+q)
\Bigr)\cdot H_k\ot1\ts,
\\
1\ot H_l\cdot Q(u)
&=
\prod_{1\leqslant p\leqslant k}^\rightarrow
\Bigl(\ts
\prod_{k<q\leqslant k+l}^\leftarrow\ts
R_{pq}(u+p+q)
\Bigr)\cdot 1\ot H_l\ts.
\endalign
$$
In particular,
$$
Q(u)\in\Fk\ot\Fl\ts[u].
$$
We will denote
$$
\phi_k\ot\phi_l\bigl(Q(u)\bigr)=\Bar Q(u).
$$
The element $\bar Q(u)$ has an inverse in
$\EndLak\ot\EndLal(u)$ as well as the element $\bar P(u)$ defined by
\(1.755).

Further, let $\ka$ and $\la$ denote the
simultaneous transpositions
in $\EndCN^{\ot(k+l)}$ with respect to the first $k$ and the last $l$
tensor factors. Due to the equalities
\(1.0) and \(3.0) the definition \(3.75) implies that
$$
\align
\ka\bigl(Q(u)\1\bigr)
&\cdot
\ka\bigl(Q(u-N)\bigr)\cdot\beta(u)=\id,
\\
\la\bigl(Q(u)\1\bigr)
&\cdot
\la\bigl(Q(u-N)\bigr)\cdot\beta(u)=\id
\endalign
$$
where $\beta(u)$ stands for the rational function
$$
\prod_{1\leqslant p\leqslant k}
\
\prod_{k<q\leqslant k+l}\ts\ts
\frac{(u+p+q)\,(u+p+q-N)}{(u+p+q-N)^2-1}\ts.
$$
Note that the elements $H_k\in\EndCN^{\ot k}$ and $H_l\in\EndCN^{\ot l}$
are invarinant with respect
to the simultaneous transpositions in all the tensor factors.
Therefore in $\EndLak\ot\EndLal(u)$ we get
$$
\align
\ka\bigl(\bar Q(u)\1\bigr)
&\cdot
\ka\bigl(\bar Q(u-N)\bigr)\cdot\beta(u)=\id,
\Tag{3.001}
\\
\la\bigl(\bar Q(u)\1\bigr)
&\cdot
\la\bigl(\bar Q(u-N)\bigr)\cdot\beta(u)=\id.
\Tag{3.02}
\endalign
$$

By \(3.2) and \(3.81) to \(3.83) we obtain in
$\EndCN^{\ot(k+l)}\!\ot\YaN\ts((u\1,v\1))$
$$
\align
&P(u-v)\ot1
\cdot
S_{1\lc k}(u,k)
\cdot
Q(-u-v)\ot1
\cdot
S_{k+1\lc k+l}(v,l)
\Tag{3.7}
\\
\qquad
&=
S_{k+1\lc k+l}(v,l)
\cdot
Q(-u-v)\ot1
\cdot
S_{1\lc k}(u,k)
\cdot
P(u-v)\ot1\ts.
\endalign
$$
Further, by using \(3.25) along with \(3.81) to \(3.83)
we obtain in $\EndCN^{\ot(k+l)}(u,v)$ the following analogue of \(3.7):
$$
\align
&P(u-v)
\cdot
Z_{1\lc k}(u,k)
\cdot
Q(-u-v)
\cdot
Z_{k+1\lc k+l}(v,l)
\Tag{3.77}
\\
&=
Z_{k+1\lc k+l}(v,l)
\cdot
Q(-u-v)
\cdot
Z_{1\lc k}(u,k)
\cdot
P(u-v)\ts.
\endalign
$$
Introduce the elements of the algebra
$\EndLak\ot\EndLal\ot\YaN[[u\1]]$
$$
\align
K(u)
&=
\phi_k\ot\phi_l\ot\id
\bigr(
\wS_{1\lc k}(u,k)
\bigl),
\\
L(u)
&=
\phi_k\ot\phi_l\ot\id
\bigl(
\wS_{k+1\lc k+l}(u,l)
\bigr)
\endalign
$$
and
$$
\align
U(u)
&=
\phi_k\bigl(Z(u+N/2\ts,k)\bigr)\ot\id\ot1\ts,
\\
V(u)
&=
\id\ot\phi_l\bigl(Z(u+N/2\ts,l)\bigr)\ot1\ts.
\endalign
\nopagebreak
$$
The equalities \(3.7) and \(3.77) then imply that
$$
\align
&
K(u)\cdot\bar Q(-u-v)\1\ot1\cdot L(v)\cdot\bar P(u-v)\ot1
\Tag{3.777}
\\
&\ =
\bar P(u-v)\ot1\cdot L(v)\cdot\bar Q(-u-v)\1\ot1\cdot K(u),
\\
&
\bar P(u-v)\ot1\cdot U(u)\cdot\bar Q(-u-v-N)\ot1\cdot V(v)=
\Tag{3.7777}
\\
&\ =
V(v)\cdot\bar Q(-u-v-N)\ot1\cdot U(u)\cdot\bar P(u-v)\ot1.
\endalign
$$
By the definition \(3.333) we have
$$
\wA_k(u)\ts\wA_l(v)=\tr\ot\id
\bigl(
K(u)\ts L(v)\cdot V(v)\ts U(u)
\bigr).
$$

Let us call two elements of
$\EndLak\ot\EndLal\ot\YaN((u\1,v\1))$
{\it equivalent} and relate them by the symbol $\sim$ if the values of
$\tr\ot\id$
on these two elements are the same. Then
$$
\align
K(u)\ts L(v)\cdot V(v)\ts U(u)
&\sim
\la\ot\id\bigl(K(u)\ts L(v)\bigr)
\cdot
\la\ot\id\bigl(V(v)\ts U(u)\bigr)
\\
&=
K(u)\cdot\la\ot\id\bigl(L(v)\bigr)
\cdot
\la\ot\id\bigl(V(v)\bigr)\cdot U(u)
\\
&=
K(u)\cdot\la\ot\id\bigl(L(v)\bigr)
\cdot
\la\bigl(\bar Q(-u-v)\1\bigr)\ot1
\\
&\hskip0.4pt\times
\la\bigl(\bar Q(-u-v-N)\bigr)\ot1
\cdot
\la\ot\id\bigl(V(v)\bigr)\cdot U(u)\cdot\beta(-u-v)
\endalign
$$
where we made use of \(3.02).
The product in the last two lines is equivalent to
$$
\align
&
\la\ot\id
\bigl(
K(u)\cdot\la\ot\id\bigl(L(v)\bigr)
\cdot
\la\bigl(\bar Q(-u-v)\1\bigr)\ot1
\bigr)
\\
\times\hskip3.5pt
&
\la\ot\id
\bigl(
\la\bigl(\bar Q(-u-v-N)\bigr)\ot1
\cdot
\la\ot\id\bigl(V(v)\bigr)\cdot U(u)
\bigr)
\cdot
\beta(-u-v)
\\
\sim\hskip3.5pt
&
K(u)\cdot\bar  Q(-u-v)\1\ot1\cdot L(v)
\\
\times\hskip3.5pt
&
V(v)\cdot\bar  Q(-u-v-N)\ot1\cdot U(u)
\cdot\beta(-u-v)
\\
=\hskip3.5pt
&
K(u)\cdot\bar Q(-u-v)\1\ot1\cdot L(v)\cdot\bar  P(u-v)\ot1
\\
\times\hskip3.5pt
&
\bar P(u-v)\1\ot1\cdot
V(v)\cdot\bar  Q(-u-v-N)\ot1\cdot U(u)
\cdot\beta(-u-v)
\\
=\hskip3.5pt
&
\bar P(u-v)\ot1\cdot L(v)\cdot\bar Q(-u-v)\1\ot1\cdot K(u)
\\
\times\hskip3.5pt
&
U(u)\cdot\bar Q(-u-v-N)\ot1\cdot V(v)\cdot\bar P(u-v)\1\ot1
\cdot\beta(-u-v)
\\
\intertext
{where we used \(3.777) and \(3.7777). The product in the last two lines
is equivalent to}
&
L(v)\cdot\bar  Q(-u-v)\1\ot1\cdot K(u)
\\
\times\hskip3.5pt
&
U(u)\cdot\bar  Q(-u-v-N)\ot1\cdot V(v)
\cdot\beta(-u-v)
\\
\sim\hskip3.5pt
&
\ka\ot\id
\bigl(
L(v)\cdot\bar Q(-u-v)\1\ot1\cdot K(u)
\bigr)
\\
\times\hskip3.5pt
&
\ka\ot\id
\bigl(
U(u)\cdot\bar Q(-u-v-N)\ot1\cdot V(v)
\bigr)\cdot
\beta(-u-v)
\\
=\hskip3.5pt
&
L(v)\cdot
\ka\ot\id
\bigl(
K(u)
\bigr)
\cdot
\ka\ot\id
\bigl(
U(u)
\bigr)
\cdot
V(v)
\sim
L(v)\ts K(u)
\cdot
U(u)\ts V(v)
\endalign
$$
where we made use of \(3.001). It now remains to take into account that
$$
\tr\ot\id
\bigl(
L(v)\ts K(u)
\cdot
U(u)\ts V(v)
\bigr)
=
\wA_l(v)\ts\wA_k(u)
\quad\square
$$
\enddemo

\proclaim{Theorem 3.5}
Suppose that the element $Z\in\EndCN$ has a simple spectrum and
$Z^{\,\prime}=-Z$. Then the subalgebra $\AaNZ$ in $\YaN$ is maximal commutative.
The coefficients at $u^{-2},u^{-4},\ldots$ of the series
$A_N(u),A_{N-2}(u),\ldots$
and
the coefficients at $u^{-1},u^{-3},\ldots$ of the series
$A_{N-1}(u),A_{N-3}(u),\ldots$
are free generators for $\AaNZ$.
\endproclaim

\nt
The proof of this theorem will be given in Section 4.
We will end up this section with making the following observation.
Consider the series $\wA(u)$ defined by \(3.333).

\proclaim{Proposition 3.6}
Suppose that $Z^{\,\prime}=\pm Z$
where the upper sign corresponds to the case $\aN=\soN$ while
the lower one corresponds to $\aN=\spN$. Then
$$
\wA_k(u)=
\tr_k\ot\id
\bigl(\ts
H_k\ot1
\cdot
\wS(u,k)
\cdot
Z_1\ldots Z_k\ot1
\ts\bigr).
$$
\endproclaim

\demo{Proof}
It can be verified directly that for $Z^{\,\prime}=\pm Z$ we have
in $\EndCN^{\ot2}[u]$ the equality
$$
Z_1\ts\tR(u)\ts Z_2\ts H_2=Z_1\ts Z_2\ts H_2\ts.
\nopagebreak
$$
By using repeatedly this equality we obtain from
\(3.007095) and \(3.333)
the required statement
\enddemos

\Section{Proof of Theorem 3.5}

\nt
We will employ arguments already used in the proof of Theorem 1.3.
Consider the ascending filtration on the algebra $\YN$
introduced in Section 2. Then by \(3.01) for the generators
of the subalgebra $\YaN$ in $\YN$ we have
$\deg S_{ij}^{(r)}=r$. Denote by $\YgaN$ and $\AgaNZ$ the
images in the graded Poisson algebra $\YgN$ of the
subalgebras $\YaN$ and $\AaNZ$ in $\YN$ respectively.
We shall prove that the subalgebra $\AgaNZ$ in the Poisson algebra
$\YgaN$ is maximal involutive.
\Par
We will denote by $s_{ij}^{(r)}$ the image in $\YgaN$
of the generator $S_{ij}^{(r)}$ of the algebra $\YaN$.
Due to the relation \(3.3) we then have
$$
s_{ij}^{(r)}=\ep_{ij}\ts s_{-j,-i}^{(r)}\cdot(-1)^r\ts;
\qquad
r=1,2,\ldots
\Tag{4.0}
$$
in $\YgaN$. Moreover, the relations
\(4.0) are defining relations for generators $s_{ij}^{(r)}$
of the commutative algebra $\YgaN$. The proof of the latter statement
is contained in [MNO\ts,\ts Section 3].
We will identify $\YgaN$ with the
symmetric algebra over the {\it twisted} polynomial current Lie algebra
$$
\bigl\{\ts F(t)\in\glN[\hskip1pt t\hskip1pt]
\ts\ts|\ts\ts\sigma\bigl(F(t)\bigr)=F(-t)\ts\bigr\}
$$
The generator $s_{ij}^{(r)}$
can be identified with the element
$$
E_{ij}\cdot t^{\,r-1}-E_{ij}^{\,\prime}\cdot(-t)^{r-1}
\Tag{4.007095}
$$
of this Lie algebra. We will set $s_{ij}^{(0)}=\de_{ij}\cdot1$.

\proclaim{Lemma 4.1}
In the Poisson algebra $\YgaN$ for any $p,q\geqslant1$ we have
$$
\align
\{\ts s_{ij}^{(p)},s_{kl}^{(q)}\ts\}
=
\ts
\sum_{r=1}^{\min(p,q)}
\bigl(\ts
s_{kj}^{(r-1)}\ts s_{il}^{(p+q-r)}-
s_{kj}^{(p+q-r)}\ts s_{il}^{(r-1)}
\ts\bigr)&
\Tag{4.1}
\\
\qquad
+\ts
\sum_{r=1}^{\min(p,q)}
\bigl(\ts
\ep_{k,-j}\ts
s_{i,-k}^{(r-1)}\ts s_{-j,l}^{(p+q-r)}-
\ep_{i,-l}\ts
s_{k,-i}^{(p+q-r)}\ts s_{-l,j}^{(r-1)}
\ts\bigr)&\cdot(-1)^{p+r-1}\ts.
\endalign
$$
\endproclaim

\demo{Proof}
By the relations \(3.22) and by the definition of the Poisson algebra
$\YgaN$
$$
\align
\{\ts s_{ij}^{(p)},s_{kl}^{(q)}\ts\}=\ts
\sum_{r=1}^p\
\bigl(\ts
s_{kj}^{(p-r)}\ts s_{il}^{(q+r-1)}-
s_{kj}^{(q+r-1)}\ts s_{il}^{(p-r)}
\ts\bigr)&
\Tag{4.01}
\\
\qquad
+\ts
\sum_{r=1}^p\
\bigl(\ts
\ep_{k,-j}\ts
s_{i,-k}^{(p-r)}\ts s_{-j,l}^{(q+r-1)}-
\ep_{i,-l}\ts
s_{k,-i}^{(q+r-1)}\ts s_{-l,j}^{(p-r)}
\ts\bigr)&\cdot(-1)^r\ts.
\endalign
$$
Here the first of the two sums coincides with the first sum in \(4.1).
The second sum in \(4.01) 
coincides with the second sum in \(4.1) by the relations \(4.0)\ \ $\square$
\enddemo

\nt
Consider the formal power series in $u\1$
$$
s_{ij}(u)=
s_{ij}^{(0)}+s_{ij}^{(1)}\ts u\1+s_{ij}^{(2)}\ts u^{-2}+\ldots
$$
with the coefficients in the algebra $\YgaN$.
The image of the series
$A_k(u)$ in the algebra $\YgaN[[u\1]]$ equals
$$
\sum_{g,\gp}\ts\ts
s_{g(1)\gp(1)}(u)
\ts\ldots\ts
s_{g(k)\gp(k)}(u)
\cdot
z_{\ts g(k+1)\gp(k+1)}\ldots\ts z_{\ts g(N)\gp(N)}
\cdot\ep\ts/\ts N!
\Tag{4.3}
$$
where $g,\gp$ run through the set of all permutations of
$-n\lc-1,1\lc n$ if $N=2n$ and the set of all permutations
$-n\lc-1,0,1\lc n$ if $N=2n+1$. Here $\ep=\sgn\ts g\cdot\sgn\ts\gp$
as well as in \(1.3) and \(2.11111111).

For each $M=1,2\ldots$ denote by $\YganN$ the ideal in the algebra
$\YgaN$ generated by all the elements $s_{ij}^{(r)}$ where $r\geqslant M$.
Note that by Lemma 4.1 we have
$$
\{\ts s_{ij}^{(p)},s_{kl}^{(q)}\ts\}\ts\in\ts\YganN\ts;
\ \quad
M=\max(p,q).
$$
Therefore each $\YganN$ is a Poisson ideal in $\YgaN$.
Now let the index $M$ be fixed.
For each $r=1\lc M$ denote by $y_{ij}^{(r)}$
the image of the generator $s_{ij}^{(r)}$in the quotient
Poisson algebra $\YgaN/\YganoN$.
In the latter algebra we have
$$
\align
\{\ts y_{ij}^{(p)},y_{kl}^{(q)}\ts\}
=\ts
\sum_{r=\max(1,p+q-M)}^{\min(p,q)}
\bigl(\ts
y_{kj}^{(r-1)}\ts y_{il}^{(p+q-r)}-
y_{kj}^{(p+q-r)}\ts y_{il}^{(r-1)}
\ts\bigr)&
\Tag{4.11}
\\
\qquad\qquad\qquad
+\ts
\sum_{r=\max(1,p+q-M)}^{\min(p,q)}
\bigl(\ts
\ep_{k,-j}\ts
y_{i,-k}^{(r-1)}\ts y_{-j,l}^{(p+q-r)}-
\ep_{i,-l}\ts
y_{k,-i}^{(p+q-r)}\ts y_{-l,j}^{(r-1)}
\ts\bigr)&\cdot(-1)^{p+r-1}
\endalign
$$
where $y_{ij}^{(0)}=1\cdot\de_{ij}\ts$.
Moreover, by \(4.0) we have the equalities
$$
y_{ij}^{(r)}=\ep_{ij}\ts y_{-j,-i}^{(r)}\cdot(-1)^r\ts;
\qquad
r=1\lc M.
$$
Let $\uMN$ be the same vector space as in \(2.2222). Consider
the subspace in $\uMN$
$$
\uaMN=
\bigl\{\ts F(t)\in\uMN
\ts\ts|\ts\ts
\sigma\bigl(F(t)\bigr)=F(-t)\ts\bigr\}.
$$
As a commutative algebra
the quotient $\YgaN/\YganoN$ can be identified with the symmetric algebra
$\operatorname{S}(\uaMN\!)$. Further, we will identify
this quotient
with the algebra $\operatorname{P}(\uaMN\!)$ of the polynomial functions
on $\uaMN$.
The generator $y_{ij}^{(r)}$ will be regarded as the coordinate
function corresponding to the vector
\(4.007095) in $\uaMN$.
\Par
Let us introduce the polynomials in $u$ of the degree $M-1$
$$
y_{ij}(u)=
y_{ij}^{(1)}\ts u^{M-1}+y_{ij}^{(2)}\ts u^{M-2}+\ts\ldots\ts+
y_{ij}^{(M)}
$$
and denote by $\hskip-1pt Y(u)\hskip-1pt$
the square matrix of order $N$ formed by these
polynomials.
For each $k=1\lc N$ denote by $a_k(u)$ the image of
$A_k(u)\in\YaN[[u\1]]$ in
$$
\YgaN/\YganoN[[u\1]]=
\operatorname{P}(\uaMN\!)[[u\1]].
$$
This image is in fact a polynomial in $u\1$
of the degree $k\hskip1pt M$, see \(4.3). We have
the Laplace expansion
$$
\det\bigl(u^M+Y(u)+Z\ts v\bigr)=
v^{N}\det Z+\sum_{k=1}^N\ts\ts
u^{kM}\ts a_k(u)\ts v^{\,N-k}
\left(\!\matrix N\\k\endmatrix\!\right)
\Tag{4.*}
$$
where $Z$ is now regarded as a square matrix of order $N$.
Since $Z^{\,\prime}=-Z$ we have
$$
a_k(-u)=a_k(u)\cdot(-1)^{N-k};
\qquad
k=1\lc N.
$$
We will write
$$
a_k(u)=a_k^{(0)}+a_k^{(1)}\ts u\1+\ts\ldots\ts+a_k^{(kM)}\ts u^{-kM}\ts.
$$
Here $a_k^{(0)}\in\CC$ for any $k=1\lc N$. Furthermore
$a_k^{(r)}=0$ unless
$
N-k+r\in2\hskip2pt\ZZ.
$
All the remaining coefficients
$a_k^{(r)}$
are in involution with respect to the Poisson bracket \(4.11)
by Proposition 3.3.

{}From now on we will be assuming that
$M=2m+1$ if $\aN=\so_{2n+1}$ or $\aN=\sp_{2n}$
and that $M=2m$ if $\aN=\so_{2n}$.
We shall prove that all the coefficients $a_k^{(r)}$ with
$$
1\leqslant r\leqslant kM\ts;
\ \quad
N-k+r\in2\hskip2pt\ZZ
\Tag{4.6}
$$
are algebraically independent and generate a maximal involutive subalgebra
in the quotient $\YgaN/\YganoN$. Lemma 2.1 will then imply that
the subalgebra $\AgaNZ$ in the Poisson algebra
$\YgaN$ is also maximal involutive. So we will then have Theorem 3.5 proved.

It can be verified directly that for any element $G\in\EndCN$
the following equality holds in the algebra $\EndCN^{\ot2}[u]$:
$$
G_1\ts\tR(u)\ts G_2^{\,\prime}
=
G_2^{\,\prime}\ts\tR(u)\ts G_1\ts;
\Tag{4.4}
$$
Now suppose that $G^{\,\prime}=G\1$ so that the element $G$
belongs to the orthogonal or symplectic group on $\CC^N$ corresponding
to the subalgebra $\aN$ in $\glN$. By the above two equalities
the defining relations \(3.2) and \(3.3) then imply that
the assignment
$$
S(u)\mapsto G\ot1\cdot S(u)\cdot G\1\ot1
$$
determines an automorphism of the algebra $\YaN$.
The image of the subalgebra
$\AaNZ$ with respect to this automorphism is
$\operatorname{B}(\glN,\sigma, G\1\hskip1pt Z\hskip1pt G)$.
Indeed, by applying \(4.4) repeatedly we obtain that the images
of the series
$A_1(u)\lc A_N(u)$ determined by the element $Z$ coincide with the
corresponding series determined by
the element $G\1\hskip1pt Z\hskip1pt G$.

Hence it suffices to prove Theorem 3.5 assuming that
$z_{ij}=z_i\cdot\de_{ij}$ where
all $z_i\in\CC$ are pairwise distinct and $z_{-i}=-z_{i}$ for any
index $i$. We shall keep to this assumption from now on
until the end of this section.

Let $\bN$ be the Borel subalgebra in $\glN$ spanned by the elements
$E_{ij}$ with $i\leqslant\nomathbreak j$.
Fix a principal nilpotent element of the Borel subalgebra in $\glN$
opposite to $\bN$
$$
E=
\cases
E_{n,n-1}-E_{1-n,-n}+\ldots+E_{21}-E_{-1,-2}+E_{1,0}-E_{0,-1}
&
\ \text{if}\ \ \aN=\so_{2n+1}\ts;
\\
E_{n,n-1}-E_{1-n,-n}+\ldots+E_{21}-E_{-1,-2}+E_{1,-1}
&
\ \text{if}\ \ \aN=\sp_{2n}\ts;
\\
E_{n,n-1}+E_{1-n,-n}+\ldots+E_{21}+E_{-1,-2}+E_{1,-1}
&
\ \text{if}\ \ \aN=\so_{2n}\ts.
\endcases
$$
\medskip
\nt
Let $\glN=\aN\oplus\xN$ be the decomposition into the eigenspaces
of the involutive automorphism $\sigma$. Then $E\in\aN$ for 
$\aN=\so_{2n+1}\ts,\sp_{2n}$ but $E\in\xN$ for 
$\aN=\so_{2n}$. Therefore
$$
E^{(M)}=E\cdot t^{\,M-1}\in\uaMN\ts
$$
by our assumption on the parity of $M$.
{We will consider the Poisson bracket \(4.11) on
$\operatorname{P}(\uaMN\!)$
in a neighbourhood of the point $E^{(M)}$\!.}

Regard the coefficients $a_k^{(r)}$ as polynomial functions
on $\uaMN$. Consider~the~affine subspace in $\uaMN$
$$
\kaN=E^{(M)}+\bigl(\ts\uaMN\cap\bN[\hskip1pt t\hskip1pt]\ts\bigr).
$$
It will be taken as an analogue of the subspace $\kN$ in $\uMN$.
For the principal nilpotent element $E\in\glN$ fixed above we have
$$
\uaMN\cap\kN=\kaN.
$$
Therefore by comparing \(2.*) and \(4.*) we
obtain the next statement directly from Lemma 2.2.

\proclaim{Lemma 4.2}
The restrictions of the functions $a_k^{(r)}$
onto the affine subspace $\kaN$
generate the whole algebra of polynomial functions on $\kaN$.
\endproclaim

\nt
Observe that the total number of the coeficients $a_k^{(r)}$ with
the indices satisfying \(4.6) coincides with
$$
\dim\ts\kaN=
\cases
(2\hskip1pt mn+m+n)(n+1)
&
\ \ \text{if}\ \ \ \aN=\so_{2n+1}\ts;
\\
(2\hskip1pt mn+m+n+1)\ts n
&
\ \ \text{if}\ \ \ \aN=\sp_{2n}\ts;
\\
(2\hskip1pt n+1)\ts m\ts n
&
\ \ \text{if}\ \ \ \aN=\so_{2n}\ts.
\endcases
$$

\proclaim{Corollary 4.3}
All the coefficients $a_k^{(r)}$ with
the indices satisfying \(4.6) are algebraically independent.
\endproclaim

\nt
To complete the proof of Theorem 3.5 we have to demonstrate
that the involutive subalgebra in $\operatorname{P}(\uaMN\!)$
generated by the coefficients $a_k^{(r)}$ is maximal.
As an argument similar to that given in the end of Section 2 shows,
this assertion follows from Lemma 4.2 and from the next lemma. Set
$$
D=\dim\ts\uaMN-\dim\ts\kaN=
\cases
(2\hskip1pt mn+m+n)\ts n
&
\ \ \text{if}\ \ \ \aN=\so_{2n+1}\ts;
\\
(2\hskip1pt mn-m+n)\ts n
&
\ \ \text{if}\ \ \ \aN=\sp_{2n}\ts;
\\
(2\hskip1pt n-1)\ts m\ts n
&
\ \ \text{if}\ \ \ \aN=\so_{2n}\ts.
\endcases
$$

\proclaim{Lemma 4.4}
There is an open neighbourhood of the point
$E^{(M)}$ in $\uaMN$
where the rank of the Poisson bracket \(4.11) equals $2D$.
\endproclaim

\demo{Proof}
Consider the alternating
bilinear form
$\alpha:\glN\times\glN\to\CC$
where $\alpha(E_{ij},E_{kl})$
is the value of \(2.4444) at the point $E$.
The rank of this bilinear form is $N^2-N$.

Any summand at the right hand side of \(4.11) vanishes at the point $E^{(M)}$
unless $p+q=M+1$ and $r=1$. We will assume that $p+q=M+1$ and that
$p,q\geqslant 1$. Then the value of the bracket
$\{\ts y_{ij}^{(p)},y_{kl}^{(q)}\ts\}$ at the point $E^{(M)}$ is
$$
\alpha\bigl(
E_{ij}-\sigma(E_{ij})\cdot(-1)^p
\ts,\ts
E_{kl}-\sigma(E_{kl})\cdot(-1)^q
\bigr).
$$

Suppose that $\aN=\so_{2n}$ so that $M=2m$.
Then $E\in\xN$ and the restriction of the bilinear form $\alpha$
onto either of
the subspaces $\aN\times\xN$ and  $\xN\times\aN$ has the rank
$$
(N^2-N)/2=n\hskip1pt(2\hskip1ptn-1).
$$
Hence there exists
an open neighbourhood of the point $E^{(M)}$ in
$\uaMN$ where
the rank of the Poisson bracket \(4.11) is not smaller than
$
2\hskip1pt m\cdot n\hskip1pt(2n-1)=2D.
$

Now suppose that $\aN=\so_{2n+1}$ or $\aN=\sp_{2n}$ so that $M=2m+1$.
Then $E\in\aN$. Moreover, $E$ is a principal nilpotent element of
the Borel subalgebra in $\aN$ opposite to $\aN\cap\bN$. Therefore
the restriction of the bilinear form $\alpha$ onto $\aN\times\aN$
has the rank $2\hskip1pt n^2$. The restriction of $\alpha$
onto $\xN\times\xN$
then has the rank
$$
(N^2-N)-2\hskip1pt n^2=2\hskip1pt n\hskip1pt (n\pm1).
$$
So there exists
an open neighbourhood of the point $E^{(M)}$ in
$\uaMN$ where
the rank of the Poisson bracket \(4.11) is not smaller than
$$
(m+1)\cdot\hskip1pt2\hskip1pt n^2+m\cdot2\hskip1pt n\hskip1pt(n\pm1)=2D.
$$

But for every $\aN$
due to Proposition 3.3 and to Corollary 4.3
in the centre of the Poisson algebra
$\operatorname{P}(\uaMN\!)$
there are algebraically independent elements
$b_N^{(r)}$ with $1\leqslant r\leqslant MN$, $r\in2\ZZ$.
So the rank of the Poisson bracket \(4.11) at any point cannot exceed
$$
2D=\dim\ts\uaMN-
\cases
2\hskip1pt m\hskip1pt n+n+m
&
\ \ \text{if}\ \ \ \aN=\so_{2n+1}\ts;
\\
2\hskip1pt m\hskip1pt n+n
&
\ \ \text{if}\ \ \ \aN=\sp_{2n}\ts;
\\
2\hskip1pt m\hskip1pt n
&
\ \ \text{if}\ \ \ \aN=\so_{2n}\ts.
\endcases
$$
The proof of Lemma 4.4 is now complete
\enddemos

\nt
Let us now consider
the commutative subalgebra $\rho\bigl(\AaNZ\bigr)$ in $\UaN$.~First
consider the involutive subalgebra in the graded algebra  
$\operatorname{S}(\aN)$ corresponding to $\rho\bigl(\AaNZ\bigr)$. 
We identify the Poisson algebras
$\operatorname{S}(\aN)$ and $\operatorname{P}(\aN)$;
the element $F_{ij}\in\aN$ is  identified with the respective
coordinate function. Denote by $y_{ij}$ this function.
Let $Y$ be
the square matrix of the order $N$ formed  all these functions.
The subalgebra  in $\operatorname{P}(\aN)$
corresponding to $\pi\bigl(\AaNZ\bigr)$
is generated by the coefficients of the polynomial in $u\ts,v$
$$
\det\bigl(u+Y+Z\ts v\bigr).
\Tag{4.999999999}
$$

This subalgebra in $\operatorname{P}(\aN)$
has been considered in [MF, Section 4].
It has been proved there that the coefficients of
\(4.999999999)
form a {\it complete set} of elements
of $\operatorname{P}(\aN)$ in involution: their
gradients at a generic point of $\aN$ span a space
of the maximal possible dimension
$$
(\dim\aN+\operatorname{rank}\aN)/2.
$$
Lemma 4.2 provides another proof of this fact for $\aN=\so_{2n+1}$
and $\aN=\sp_{2n}$. In these two cases it also shows that the
coefficients of \(4.999999999)
generate a maximal involutive subalgebra in $\operatorname{P}(\aN)$.
Indeed, when $\aN=\so_{2n+1}$ or $\aN=\sp_{2n}$ we can set $M=1$.
But when $\aN=\so_{2n}$ by our assumption $M$ has to be even.
So we have to consider the latter case separately.

Suppose that $N=2n$ and $\aN=\so_{2n}$. Fix in
$\so_{2n}$
the Borel subalgebra $\frak b_{2n}\cap\so_{2n}$
and the principal nilpotent element
of the opposite Borel subalgebra
$$
E=E_{n,n-1}-E_{1-n,-n}+\ldots+E_{21}-E_{-1,-2}+E_{2,-1}-E_{1,-2}\ts.
$$
Consider the affine subspace
$\frak{s}_{2n}=E+\frak b_{2n}\cap\so_{2n}$ in $\so_{2n}$.
To prove that the
coefficients of \(4.999999999) again
generate a maximal involutive subalgebra in $\operatorname{P}(\so_{2n})$
it suffices to establish the following lemma.

\proclaim{Lemma 4.5}
The restrictions
of the coefficients of the polynomial \(4.999999999)
onto the affine subspace $\frak{s}_{2n}$
generate the whole algebra of polynomial functions on $\frak{s}_{2n}$.
\endproclaim

\demo{Proof}
We will employ arguments already used in the proof of Lemma 2.2.
Denote
the polynomial in $u,v$ obtained by restricting the coefficients of
\(4.999999999) onto $\frak{s}_{2n}$ by $F(u,v)$. Let the indices $i,j$
run through the set $\{1\lc n\}$.
Take the functions $y_{ij}$ with $i\leqslant j$
and $y_{-i,j}$ as coordinates on $\frak{s}_{2n}$.
Put $\bar y_{ij}=y_{ij}-\de_{1i}\ts y_{-1,j}$.
We will prove consecutively for $d=0,1\lc 2n-2$ that here
all $\bar y_{ij}$ with $j-i=d$ and all $y_{-i,j}$ with $i+j=d+2$
are polynomials in the coefficients of $F(u,v)$.

Consider the terms of $F(u,v)$ with the total degree $2n-d-1$ in $u,v$.
If $d=0\lc n-1$ their sum has the form
$$
f(u,v)+(-1)^{\,n+d}\cdot\ts
\sum_{i=1}^{n-d}\ts\ts
\bar y_{i,i+d}
\ts\ts g_i(u,v)
-
(-1)^{\,n+d}\cdot\ts\sum_{i=1}^{d+1}\ts\ts
y_{-i,d-i+2}
\ts\ts h_i(u,v)
$$
where
$$
h_i(u,v)=
2u
\ts\ts\cdot
\prod_{j=i+1}^n\ts (-\,u+z_j\ts v)
\ts\ts\cdot\!
\prod_{j=d-i+2}^n\ts (u+z_j\ts v)
\Tag{4.**}
$$
while $g_i(u,v)$ equals
$$
\prod_{j\neq i\lc i+d}\ts (-\,u+z_j\ts v)
\ts\cdot\ts
\prod_{j=1}^n\ts\ts (u+z_j\ts v)
\ts+\!
\prod_{j\neq i\lc i+d}\ts (u+z_j\ts v)
\ts\cdot\ts
\prod_{j=1}^n\ts\ts (-\,u+z_j\ts v).
$$
Here the coefficients of the polynomial $f(u,v)$
depend only on $\bar y_{ij}$ with $j-i<d$ and on $y_{-i,j}$ with $i+j<d+2$
along with $z_1\lc z_n$.
But all the $n+1$
polynomials
$g_i(u,v)$ and $h_i(u,v)$ are linearly independent.
To prove this consider their values at
$$
(u,v)=(z_1\ts,1)\lc(z_n\ts,1),\ts(0\ts,1).
$$
Since
all the numbers $z_1\ts,-z_1\ts\lc z_n\ts,-z_n$ are pairwise distinct
we have
$$
\align
h_i(z_i\ts,1)&\neq0\ts,
\qquad
1\leqslant i\leqslant d+1\ts;
\\
h_j(z_i\ts,1)&=0\ts,
\qquad
1\leqslant j<i\leqslant d+1\ts.
\\
\intertext{Moreover, we have}
h_j(0\ts,1)&=0\ts,
\qquad
1\leqslant j\leqslant d+1\ts;
\\
h_j(z_i\ts,1)&=0\ts,
\qquad
1\leqslant j\leqslant d+1<i\leqslant n\ts.
\\
\intertext{Furthermore, we have $g_1(0,1)\neq0$ and}
g_i(z_{d+i}\ts,1)&\neq0\ts,
\qquad
2\leqslant i\leqslant n-d\ts;
\\
g_j(z_{d+i}\ts,1)&=0\ts,
\qquad
1\leqslant j<i\leqslant n-d\ts.
\endalign
$$

If $d=n\lc 2n-2$ then the sum of the terms of $F(u,v)$
with the total degree $2n-d-1$ has the form
$$
f(u,v)-(-1)^{\,n+d}\cdot\ts
\sum_{i=d-n+2}^{n}\ts\ts
y_{-i,d-i+2}
\ts\ts h_i(u,v)
$$
where $h_i(u,v)$ is determined by \(4.**) while
the coefficients of the polynomial $f(u,v)$
depend only on $\bar y_{-i,j}$ with $i+j<d+2$ and on $\bar y_{ij}$
along with $z_1\lc z_n$.
Here all the $2n-d-1$ polynomials $h_i(u,v)$ are again linearly
independent. Indeed,
$$
\align
h_i(z_i\ts,1)&\neq0\ts,
\qquad
d-n+2\leqslant i\leqslant n\ts;
\\
h_j(z_i\ts,1)&=0\ts,
\qquad
d-n+2\leqslant j<i\leqslant n
\qquad\square
\endalign
$$

\nt
Thus for
$\aN=\so_{2n+1}\ts,\ts\sp_{2n}\ts,\ts\so_{2n}$ we have proved that
the involutive subalgebra~in $\operatorname{P}(\aN)$ 
corresponding to the commutative subalgebra $\rho\bigl(\AaNZ\bigr)$ of
$\UaN$,
is maximal.
So for each of the above classical Lie algebras
we get the next theorem.

\proclaim{Theorem 4.6}
The subalgebra $\rho\bigl(\AaNZ\bigr)$ in $\UaN$
is {maximal commutative.}
\endproclaim

\nt
This theorem is the analogue of Proposition 2.6 for the above classical Lie
algebras.

\section{Acknowledgements}

\nt
The first author
was supported by the University of Wales
Research Fellowship.
The second author
was supported by
Russian Foundation for Basic Research
under Grant 95-01-00814.

\section{References}

\medskip

\itemitem{[C1]}
{I. V. Cherednik},
{\it On the properties of factorized $S$ matrices in elliptic functions},
{Sov.\,J.\,Nucl.\,Phys.}
{\bf 36}
(1982),
320--324.

\itemitem{[C2]}
{I. V. Cherednik},
{\it Factorized particles on the half-line and root systems},
{Theor. Math.\,Phys.}
{\bf 61}
(1984),
35--44.

\itemitem{[C3]}
{I. V. Cherednik},
{\it A new interpretation of Gelfand-Zetlin bases},
{Duke Math.\,J.}
{\bf 54}
(1987),
563--577.

\itemitem{[D1]}
{V. G. Drinfeld},
{\it Hopf algebras and the quantum Yang--Baxter equation},
{\text{Soviet} Math.\,Dokl.}
{\bf 32}
(1985),
254--258.

\itemitem{[D2]}
{V. G. Drinfeld},
{\it A new realization of Yangians and quantized affine algebras},
{Soviet Math. Dokl.}
{\bf 36}
(1988),
212--216.

\itemitem{[H]}
{R. Howe},
{\it Some Highly Symmetrical Dynamical Systems},
preprint.

\itemitem{[JM]}
{M. Jimbo and T. Miwa},
{\it Algebraic Analysis of Solvable Lattice Models},
Regional Conference Series in Math.
{\bf 85},
AMS, Providence RI, 1995.

\itemitem{[K1]}
{B. Kostant},
{\it Lie group representations on polynomial rings},
{Amer.\,J.\,Math.}
{\bf 85}
(1963),
327--404.

\itemitem{[K2]}
{B. Kostant},
{\it The solution to a generalized Toda lattice and representation theory,}
{Adv.\,Math.}
{\bf 34}
(1979),
195--338.

\itemitem{[KBI]}
{V. E. Korepin, N. M. Bogoliubov and A. G. Izergin},
{\it Quantum Inverse Scattering Method and Correlation Functions},
Cambridge University Press, Cambridge, 1993.

\itemitem{[KR]}
{A. N. Kirillov {\rm and} N. Yu. Reshetikhin},
{\it Yangians, Bethe ansatz and combinatorics},
{Lett.\,Math.\,Phys.}
{\bf 12}
(1986),
199--208.

\itemitem{[KS]}
{P. P. Kulish {and} E. K. Sklyanin},
{\it Quantum spectral transform method: recent developments},
in \lq\ts Integrable Quantum Field Theories\ts\rq,
{Lecture Notes in Phys.}
{\bf 151},
Springer,
Berlin--Heidelberg,
1982,
pp. 61--119.

\itemitem{[M]}
{A. Molev},
{\it Sklyanin determinant, Laplace operators
and characteristic identities for classical Lie algebras},
{J.\,Math.\,Phys.}
{\bf 36}
(1995),
923--943.

\itemitem{[MF]}
{A. S. Mishchenko and A. T. Fomenko},
{\it Euler equations on finite-dimensional Lie groups},
{Izv.\,AN SSSR Ser.\,Math.}
{\bf 42}
(1978),
396--415.

\itemitem{[MNO]}
{A. Molev, M. Nazarov and G. Olshanski\u\i},
{\it Yangians and classical Lie algebras},
{Preprint CMA}
{\bf 53},
Austral.\,Nat.\,Univ.,
Canberra,
1993;
{\tt hep-th/9409025}.

\itemitem{[NS]}
{M. Noumi and T. Sugitani},
{\it Quantum symmetric spaces and related $q$-ortho\-gonal polynomials},
{Preprint UTMS}
{\bf},
Univ.\,of Tokyo,
1994.

\itemitem{[NT]}
{M. Nazarov and V. Tarasov},
{\it Representations of Yangians with Gelfand-Zetlin bases},
{Preprint MRRS}
{\bf 148},
Univ.\,of Wales, Swansea,
1994;
{\tt q-alg/9502008}.

\itemitem{[O1]}
{G. I. Olshanski\u\i},
{\it Representations of infinite-dimensional classical
groups, limits of enveloping algebras, and Yangians},
in \lq\ts Topics in Representation Theory\ts\rq,
{Advances in Soviet Math.}
{\bf 2},
AMS,
Providence,
1991,
pp. 1--66.

\itemitem{[O2]}
{G. I. Olshanski\u\i},
{\it Twisted Yangians and infinite-dimensional classical Lie algebras},
in
\lq\ts Quantum Groups\ts\rq,
{Lecture Notes in Math.}
{\bf 1510},
Springer,
Berlin--Heidelberg,
1992,
pp. 103--120.

\itemitem{[RS]}
{A. G. Reiman and M. A. Semenov-Tian-Shansky},
{\it Reduction of Hamiltonian systems, affine Lie algebras and Lax equations
II,}
{Invent.\,Math.}
{\bf 63}
(1981),
423--432.

\itemitem{[RT]}
{M. Ra{\"\i}s and P. Tauvel},
{\it Indice et polyn\^omes invariants pour certaines alg\`ebres de Lie},
{J.\,Reine Angew.\,Math.}
{\bf 425}
(1992),
123--140.
 
\itemitem{[S]}
{E. K. Sklyanin},
{\it Boundary conditions for integrable quantum systems}, {J.\,Phys.}
{\bf A21}
(1988),
2375--2389.

\itemitem{[V]}
{E. B. Vinberg},
{\it On certain commutative subalgebras of a universal enveloping algebra},
{Izv.\,AN SSSR Ser.\,Math.}
{\bf 54}
(1990),
3--25.

\bye